\documentclass[11pt]{article}
\pdfoutput=1                       

\usepackage[utf8]{inputenc}        
\usepackage[english]{babel}        
\usepackage[hidelinks]{hyperref}   
\usepackage{fullpage}              
\usepackage{parskip}               
\usepackage[shortlabels]{enumitem} 
\usepackage{cite}                  
\usepackage[margin=2em,labelfont=bf,size=small]{caption}
\usepackage{amsmath,amsthm,amsfonts,amssymb}    
\usepackage[capitalise,nameinlink]{cleveref}	
\usepackage{booktabs}              
\usepackage[bottom]{footmisc}

\usepackage{graphicx}
\usepackage{xcolor}
\graphicspath{{./graphics/}}

\usepackage{bm}  

\usepackage{colonequals}                      
\usepackage{bookmark}                         

\DeclareMathOperator*{\argmin}{\arg\min}
\DeclareMathOperator{\E}{\mathbb{E}}
\DeclareMathOperator{\trace}{\mathrm{tr}}
\DeclareMathOperator{\diag}{\mathrm{diag}}

\newcommand{\bmat}[1]{\begin{bmatrix}#1\end{bmatrix}}                        
\newcommand{\sbmat}[1]{\left[\begin{smallmatrix}#1\end{smallmatrix}\right]}  
\newcommand{\R}{\mathbb{R}}					       
\newcommand{\tp}{\mathsf{T}}                       
\newcommand{\defeq}{\colonequals}                  
\newcommand{\e}{\mathsf{e}}                        
\newcommand{\F}{\mathcal{F}}                       
\newcommand{\df}{\nabla\! f}                       
\newcommand{\norm}[1]{\lVert{#1}\rVert}            
\newcommand{\normm}[1]{\bigl\lVert{#1}\bigr\rVert} 
\newcommand{\alg}{\mathcal{A}}
\newcommand{\col}{\text{col}}
\newcommand{\y}{\bm{y}}
\renewcommand{\u}{\bm{u}}
\newcommand{\bxi}{\bm{\xi}}
\newcommand{\x}{\bm{x}}
\newcommand{\f}{\bm{f}}
\newcommand{\A}{\bm{A}}
\newcommand{\B}{\bm{B}}
\newcommand{\C}{\bm{C}}
\newcommand{\D}{\bm{D}}
\renewcommand{\H}{\bm{H}}

\newcommand{\X}{\bm{X}}
\newcommand{\Y}{\bm{Y}}
\newcommand{\U}{\bm{U}}

\newcommand{\rate}{\textup{\textsc{rate}}}
\newcommand{\sensitivity}{\textup{\textsc{sens}}}

\renewcommand{\epsilon}{\varepsilon}
\newcommand{\squeezemat}{\addtolength{\arraycolsep}{-2pt}}

\makeatletter
\def\thm@space@setup{%
  \thm@preskip=6pt plus 2pt minus 1pt
}
\makeatother

\newtheorem{theorem}{Theorem}[section]
\newtheorem{lemma}[theorem]{Lemma}
\newtheorem{corollary}[theorem]{Corollary}
\newtheorem{proposition}[theorem]{Proposition}

\newtheorem{remark}{Remark}

\newtheorem{assumption}{Assumption}

\def\qed{\hfill~\rule[0pt]{5pt}{5pt}\par\medskip}
\renewenvironment{proof}{{\noindent\bf Proof.}}{\qed}

\newcommand{\parag}[1]{\paragraph*{#1.}}

\begin{document}

\title{The Speed--Robustness Trade-Off for First-Order Methods\\
with Additive Gradient Noise}

\author{Bryan Van Scoy\thanks{B. Van Scoy is with the Department of Electrical and Computer Engineering, Miami University,\newline
Oxford, OH, USA. Email: \texttt{bvanscoy@miamioh.edu}}
\and
Laurent Lessard\thanks{L. Lessard is with the Department of Mechanical and Industrial Engineering, Northeastern University,\newline
Boston, MA, USA. Email: \texttt{l.lessard@northeastern.edu}}
}
\date{}

\maketitle

\begin{abstract}
  We study the trade-off between convergence rate and sensitivity to stochastic additive gradient noise for first-order optimization methods. Ordinary Gradient Descent (GD) can be made fast-and-sensitive or slow-and-robust by increasing or decreasing the stepsize, respectively.	However, it is not clear how such a trade-off can be navigated when working with accelerated methods such as Polyak's Heavy Ball (HB) or Nesterov's Fast Gradient (FG) methods.	We consider two classes of functions: (1) strongly convex quadratics and (2) smooth strongly convex functions. For each function class, we present a tractable way to compute the convergence rate and sensitivity to additive gradient noise for a broad family of first-order methods, and we present algorithm designs that trade off these competing performance metrics. Each design consists of a simple analytic update rule with two states of memory, similar to HB and FG. Moreover, each design has a scalar tuning parameter that explicitly trades off convergence rate and sensitivity to additive gradient noise. We numerically validate the performance of our designs by comparing their convergence rate and sensitivity to those of many other algorithms, and through simulations on Nesterov's ``bad function''.
\end{abstract}

\section{Introduction}\label{sec:intro}

We consider the problem of designing robust first-order methods for unconstrained minimization. Given a continuously differentiable function $f : \R^d\to\R$, consider solving the optimization problem
\begin{equation}\label{eq:opt}
  x^\star \in \argmin_{x\in\R^d} f(x),
\end{equation}
where the algorithm only has access to gradient measurements corrupted by \textit{additive stochastic noise}.\footnote{%
Preliminary versions of portions of this work appeared in the conference proceedings \cite{alift_cdc,tut-lyaplift_cdc}.} Specifically, the algorithm can sample the oracle $g(x) \defeq \df(x) + w$, where $w$ is zero-mean and independent across queries. This form of additive noise arises in various applications. For example, (i) to protect sensitive data, optimization algorithms may intentionally perturb the gradient by Gaussian noise in order to obtain differential privacy~\cite{bassily2014}; (ii) for some engineering systems, the gradient can only be obtained through noisy measurements~\cite{birand2013}; and (iii) in risk minimization in the context of learning algorithms, the objective is to minimize the expectation of the loss function over the population distribution~\cite{lan2012,ghadimi-lan1,kakade_accelerated_stochastic}.

\begin{figure}[ht]
	\centering
	\includegraphics[width=\linewidth]{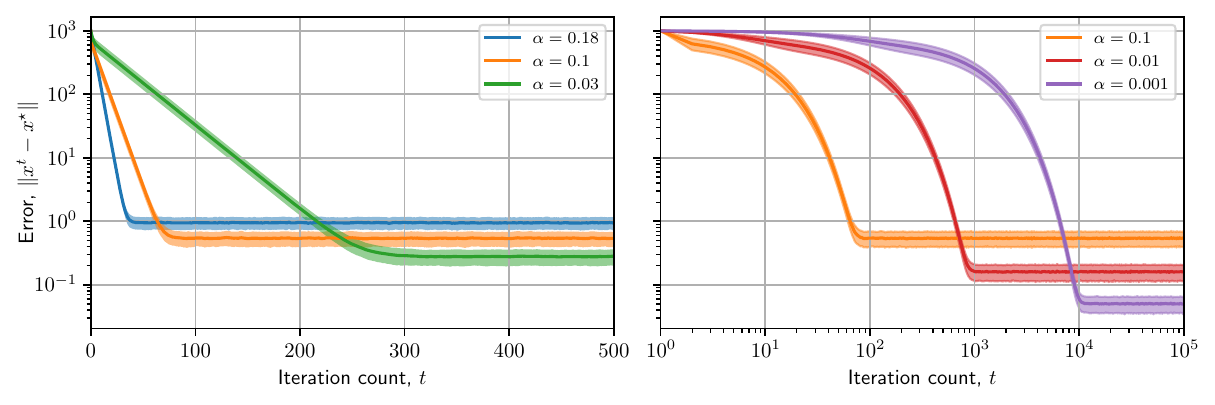}
    
\caption{Trade-off between convergence rate and steady-state error (sensitivity to noise). Three different tunings of Gradient Descent (GD) with additive gradient noise are applied to random strongly convex quadratic functions on $\R^{10}$. Half of the Hessian eigenvalues are at $m=1$, the other half at $L=10$. The initialization is $x^0 = 1000\,\e_1$. Gradient noise is normally distributed $\mathcal{N}(0,I)$ and i.i.d.\ across iterations. The plot shows mean and $\pm 1$ standard deviation of the error $\norm{x^t-x^\star}$ for 1000 sample trajectories.
On the right panel, iterations are plotted on a log scale.}
\label{fig:demo_rand_trajectories}
\end{figure}

Many iterative algorithms have been proposed to solve this problem, and most have tunable parameters. For example, Gradient Descent (GD) uses the update
\begin{align}\label{eq:sgd}
	\text{Gradient Descent (GD):} & \qquad x^{t+1} = x^t - \alpha\, g(x^t),
\end{align}
where $t$ is the iteration index and the stepsize $\alpha$ is a tunable parameter. \cref{fig:demo_rand_trajectories} illustrates how the error $\norm{x^t-x^\star}$ evolves under GD applied to strongly convex quadratic functions for different fixed $\alpha$.
The convergence of the error is characterized by an initial transient phase followed by a stationary phase. In the transient phase, the gradient dominates the noise, and the error converges at a linear rate. When the gradient is small enough that the noise becomes significant, the average error of the iterates converges to a constant value.

This two-phase behavior is typical of stochastic methods.\footnote{In the literature, stepsize is also known as \emph{learning rate}. The transient phase is also known as the \emph{search} or \emph{burn-in} phase. The stationary phase is also known as the \emph{convergence} or \emph{steady-state} phase.\looseness=-1}
The fundamental trade-off is that faster initial convergence comes at the cost of larger steady-state~error.

Gradient Descent is easy to interpret and tune: the stepsize directly mediates the trade-off between convergence rate and sensitivity. Unfortunately, GD is generally slow to converge, and alternative methods can provide faster convergence rates. Two such methods are Polyak's Heavy Ball~\cite{polyak_book} and Nesterov's Fast Gradient~\cite{nesterov_book}, which use the updates
\begin{align}
	\text{Heavy Ball (HB):} &\quad x^{t+1} = x^t - \alpha\,g(x^t) + \beta\,(x^t - x^{t-1}) \label{eq:hb}, \\
  \text{Fast Gradient (FG):} &\quad x^{t+1} = x^t - \alpha\,g\bigl(x^t + \beta\,(x^t - x^{t-1})\bigr) + \beta\,(x^t - x^{t-1}).\label{eq:fg}
\end{align}
When using a noisy gradient oracle, these methods exhibit a stationary phase similar to that of GD in \cref{fig:demo_rand_trajectories}.
A trade-off between convergence rate and sensitivity to noise must also exist for HB and FG, but there are now two parameters to tune, so it is unclear how they should be modified to mediate this trade-off.

The goal of the present work is to study the trade-off between rate and sensitivity for first-order algorithms, and to design algorithms that trade off these competing performance metrics.
We consider two well-studied classes of functions $f:\R^d\to\R$ characterized by scalar parameters $m$ and $L$ that satisfy $0 < m \leq L < \infty$.

\begin{itemize}[itemsep=1mm,topsep=1mm]
	\item \textbf{Strongly convex quadratics ($Q_{m,L}$).}
    Functions such that, for some $f^\star \in \R$, $y^\star \in \R^d$, and $Q=Q^\tp\in\R^{d\times d}$ with all eigenvalues in the interval $[m,L]$, the function has the form $f(y) = \tfrac{1}{2} (y-y^\star)^\tp Q (y-y^\star)+f^\star$.
	\item \textbf{Smooth strongly convex functions ($F_{m,L}$).}
	Continuously differentiable functions for which $f(y) \geq f(x)+\df(x)^\tp (y-x) + \tfrac{m}{2}\norm{y-x}^2$ (strongly convex) 
    and $\norm{\df(x)-\df(y)} \leq L\norm{x-y}$ (Lipschitz gradients) for all $x,y\in\R^d$.
\end{itemize}
These classes of functions are nested in the sense that
$Q_{m,L} \subseteq F_{m,L}$ with equality occurring when $m=L$. Moreover, all functions in these classes have
a unique optimal point $y^\star = \argmin_{y\in\R^d} f(y)$ satisfying $\df(y^\star)=0$.

For each function class, we design algorithms in the three-parameter family:
\begin{equation}\label{eq:algform0}
	x^{t+1} = x^t - \alpha\,g\bigl( x^t + \eta\,(x^t-x^{t-1})\bigr) + \beta\,(x^t - x^{t-1}),
\end{equation}
where $\alpha$, $\beta$, $\eta$ are scalar parameters. This parameterization was first introduced in~\cite{lessard16} and is further discussed in \cref{sec:parameterization}. Note that GD, HB, FG are special cases of the three-parameter family~\eqref{eq:algform0} that use $\beta=\eta=0$, $\eta=0$, and $\eta=\beta$, respectively.

\subsection{Main contributions}

Our three main contributions are as follows.
\begin{enumerate}
	\item For the function classes $Q_{m,L}$ and $F_{m,L}$, we present a method for efficiently bounding the worst-case convergence rate and sensitivity to additive gradient noise for a wide class of algorithms. The computational effort required to find the rate and sensitivity bounds for a given algorithm is independent of problem dimension $d$ and takes fractions of a second on a laptop.\looseness=-1
	\item We present two new algorithms, RHB and RAM, for the function classes $Q_{m,L}$ and $F_{m,L}$, respectively. Each algorithm has the form~\eqref{eq:algform0}, where $(\alpha,\beta,\eta)$ are explicit algebraic functions of $m$ and $L$ and a scalar tuning parameter $r$ that  directly trades off convergence rate and sensitivity to gradient noise.
	\item We demonstrate through empirical studies that our algorithm designs effectively trade-off convergence rate and sensitivity to gradient noise. We use a brute-force approach to compare our algorithms against a dense sampling of algorithms of the form \eqref{eq:algform0}. We also show that our designs compare favorably to (i) popular algorithms such as nonlinear conjugate gradient and quasi-Newton methods, and (ii) existing designs that use more parameters.
\end{enumerate}

\parag{Paper organization} We first provide additional background on the problem and describe our analysis and design methodology. In \cref{sec:defs}, we describe the problem setting and performance measures we will use. \Cref{sec:Q,sec:F} treat $Q_{m,L}$ and $F_{m,L}$, respectively. For each, we present a computationally tractable approach to computing convergence rate and noise sensitivity, and a novel algorithm design. In \cref{sec:numerical}, we present empirical studies that support the effectiveness of our designs and compare them to existing algorithms. We conclude in \cref{sec:conclusion}.

\subsection{Literature review}\label{sec:background}

\parag{Complexity bounds}

The noise-free setting is well-studied, and algorithms have been discovered that achieve optimal worst-case iteration complexity for a variety of different function classes. 
We discuss these results in \cref{sec:algodesign_QmL,sec:algodesign_FmL}, where we present our algorithm designs for $Q_{m,L}$ and $F_{m,L}$, respectively.

In the additive gradient noise setting, there are fundamental lower bounds on the asymptotic convergence rate of iterative methods. No matter what iterative scheme is used, $\mathbb{E}\norm{x^t-x^\star}^2$ cannot decay asymptotically to zero faster than $1/t$. Roughly, this is because the rate at which error can decay is limited by the concentration properties of the gradient noise \cite{ghadimi-lan1,ghadimi-lan2}. This optimal asymptotic rate is achieved by Gradient Descent with a \emph{diminishing stepsize} that decays as $1/t$ \cite[Thm.~4.7]{bottou2018} (GDDS).
However, GDDS can suffer from poor finite-time performance: iterates may require many steps to reach the asymptotic regime. In practice, algorithms are run for a finite horizon, and if the transient phase persists, the resulting error can remain large. Hence, asymptotic optimality does not imply finite-time efficiency.

\parag{Stochastic approximation}

In the additive gradient noise setting, a variety of techniques have been proposed to improve transient performance. For the case of least squares regression, a dramatic speedup can be achieved via careful manipulation of the stepsize or by using acceleration~\cite{kakade_accelerated_stochastic,kakade_step_decay}. Kulunchakov and Mairal~\cite{kulunchakov2019} also show that any algorithm that converges linearly for smooth and strongly convex objectives in the noiseless setting can be converted into an algorithm that converges optimally (up to log factors) when additive noise is included.

Ghadimi and Lan~\cite{ghadimi-lan1} propose the AC-SA algorithm, which is a modification of Nesterov's Fast Gradient method that has near-optimal iteration complexity in the smooth and strongly convex setting.
Besides characterizing the convergence of the average iterates, they also estimate the performance of the algorithm on a single problem instance.
Building on this work, Cohen et al.~\cite{cohen2018} propose the algorithm AGD+, a ``lazy'' (dual averaging) counterpart of AC-SA for the smooth and convex setting, along with its extension $\mu$AGD+ for the strongly convex setting. In the strongly convex setting, $\mu$AGD+ improves on the convergence bound of the AC-SA algorithm, but does not achieve the lower complexity bound. In particular, applying \cite[Corollary B.5]{cohen2018} with $\gamma_i = c\sqrt{m/L}$ yields the bound
\begin{equation}\label{eq:AGD}
    \frac{m}{2} \E\|y^t-y^\star\|^2
    \;\leq\; \left(1-c\sqrt{\frac{m}{L}}\right)^t \frac{L-m}{2} \|x^0-x^\star\|^2 + \frac{\sigma^2 c^2}{\sqrt{mL}},
\end{equation}
where we used the definition of strong convexity to bound the distance to optimality in terms of the optimality gap. The parameter $c\in(0,\sqrt{L/m}]$ can be chosen to trade off rate and sensitivity for this algorithm; however, we show that the resulting trade-off is strictly suboptimal (see \cref{fig:FmL_cloud}).

None of the above algorithms achieve the optimal iteration complexity in the smooth strongly convex setting. One way to achieve the optimal complexity, however, is to exploit the rapid convergence of the transient phase by using a piecewise constant stepsize, effectively dividing the convergence into \emph{epochs} or \emph{stages}. This can be done on a predetermined schedule, or by using a statistical test to detect the phase transition which then triggers the parameter change \cite{toulis_convergence_diagnostic,cohen2018}. In Ghadimi and Lan's companion paper~\cite{ghadimi-lan2}, for instance, they propose a multi-stage version of AC-SA that achieves the optimal expected rate of convergence (up to constants) in this setting. Aybat et al.~\cite{M-ASG} build on Nesterov's Accelerated Stochastic Gradient method to construct the Multi-Stage Accelerated Stochastic Gradient method (M-ASG). This algorithm simultaneously achieves optimal iteration complexity (again, up to constants) in both deterministic and stochastic cases without knowledge of the noise parameters.\looseness=-1

The aforementioned multi-stage algorithms tune the single-stage algorithm stepsizes and the restart schedule to achieve desirable transient and asymptotic convergence, and they achieve the optimal complexity up to constants.
In the present work work, we focus on analyzing and designing single-stage algorithms with fixed stepsizes that can be tuned via a single scalar parameter to explicitly trade-off convergence rate and sensitivity to noise.
In other words, our algorithms can easily be adjusted to be made faster (and more sensitive to noise), or more robust to noise (and slower). When tuned to be as fast as possible, our algorithm designs recover algorithms that are known to converge with the optimal linear rate in the noise-free setting for the respective function classes (see \cref{sec:numerical}). Our algorithm designs may be used to converge at a linear rate to a predetermined noise level, or the trade-off parameter may be adjusted over time to construct multi-stage variants\footnote{While we do not propose a restart mechanism or time-varying parameter schedule in the present work, we illustrate such an approach in \cref{sec:simulation_worst-case} through a hand-tuned schedule.}.

\parag{Other noise models}

While we restrict our attention to additive stochastic gradient noise in the present work, other inexact oracles have been studied in the literature.
When the gradient noise is due to sampling a finite-sum objective function, viable strategies include incremental gradient or variance reduction methods such as SVRG \cite{svrg}, SAGA \cite{saga}, and many others \cite{wang2013,spider,katyusha,snvrg,sarah}.

A prominent alternative model is to assume \emph{deterministic}\footnote{Also known as \emph{worst-case} or \emph{adversarial} noise.} bounded noise, which may be additive or multiplicative.
For instance, Devolder et al.~\cite{dgn} propose an inexact deterministic first-order oracle and show that, under this oracle, acceleration necessarily results in an accumulation of gradient errors in the unconstrained smooth (not strongly convex) case.
For this setting, the same authors propose the Intermediate Gradient Method, a family of first-order methods that can be tuned to trade off convergence rate and sensitivity to gradient errors~\cite{IGM}.
The same authors also extended this oracle to analyze inexact first-order methods in the strongly convex case~\cite{devolder-strongly-convex}; there is still a trade-off between rate and sensitivity, but errors no longer accumulate and can converge linearly to within a constant ball about the optimal solution. Multiplicative deterministic noise has been studied in~\cite{lessard16}, for which the Robust Momentum (RM)~\cite{rmm} was designed to trade off convergence rate and sensitivity. Stochastic Gradient Descent has also been analyzed under a hybrid additive and multiplicative deterministic oracle~\cite{sgdn}.

\subsection{Methodology}
\label{sec:design}

We now provide an overview of our analysis and design framework and point to related techniques.
Our method for bounding the convergence rate $r$ and sensitivity $\gamma$ for the class $F_{m,L}$ relies on solving a small linear matrix inequality (LMI). This idea builds upon several related works.

The \textit{performance estimation problem} (PEP) framework~\cite{drori-teboulle,taylor2017smooth} uses LMIs to directly search for a problem instance that achieves the worst-case performance. The PEP framework has been successfully applied to numerous algorithms, such as the proximal gradient method~\cite{taylor2018exact}, operator splitting methods~\cite{ryu2020}, and gradient descent using an exact line search with noisy search directions~\cite{deklerk2017}. Of particular interest, \cite{taylor2019} uses PEP to analyze iterative algorithms under various noise models when the objective is smooth and convex. This approach searches for a time-varying potential function whose existence certifies a performance bound. To obtain closed-form bounds, increasingly large LMIs are solved and their numerical solutions guide the construction of explicit expressions for the potential function parameters as functions of the iteration index.
Alternatively, small PEPs can be formulated and solved to search for potential functions with fixed parameters as in~\cite{deKlerk2020}.

Viewing algorithms as discrete-time dynamical systems, Integral Quadratic Constraints (IQCs) from control theory~\cite{lessard16,scherer,fazlyab} may be used to search for worst-case performance guarantees. This approach also leads to LMIs that characterize \textit{asymptotic} performance, although the ensuing performance bounds are not tight in general.
  
Finally, LMIs may be used to directly search for a \textit{Lyapunov function}, which is a generalized notion of ``energy'' stored in the system. If a fraction of the energy dissipates at each iteration, this is akin to proving convergence at a specified rate \cite{diss_ICML,taylor2018lyapunov,lessardCSM}. Similar to IQCs, Lyapunov functions provide asymptotic (albeit more interpretable) performance guarantees.

In the present work, we adopt a Lyapunov approach most similar to~\cite{taylor2018lyapunov}, but we generalize it to include both convergence rate and sensitivity to noise. We also explain in \cref{sec:comparison} how the Lyapunov approach is related to IQCs.

Given an LMI that establishes the performance (convergence rate or sensitivity to noise) of a given first-order method, the next step is to \textit{design} algorithms that exploit this trade-off. Several approaches have been proposed: i) One can parameterize a family of candidate algorithms and search over parameters to effectively trade off convergence rate and sensitivity to noise.
Such a problem is typically non-convex, so one must resort to exhaustive search \cite{lessard16} or nonlinear numerical solvers that find local optima.
ii) Using convex relaxations or other heuristics such as coordinate descent, the algorithm design problem can be solved approximately. While this approach may lead to conservative designs, it has the benefit of being automated, flexible, and amenable to efficient convex solvers~\cite{scherer}.
iii) In certain settings, the controller parameters can be eliminated from the LMI entirely, yielding bounds that hold for all algorithms. This approach has been used to show that the Triple Momentum (TM) method achieves the optimal worst-case rate over the class $F_{m,L}$ \cite{algosyn_ACC,scherer2021convex}.

As in the first approach, we use the parameterized family of algorithms \eqref{eq:algform0} and search over the parameters. However, our approach to design is algebraic rather than numeric. We find analytic solutions to the non-convex semidefinite programs that arise when the algorithm parameters are treated as decision variables. This approach enables us to find explicit analytic expressions for our algorithm parameters.
Nevertheless, we still make use of numerical approaches in order to validate our choice of parameterization and the efficacy of our designs, and to compare the performance of our designs to that of existing algorithms; see \cref{sec:numerical}.

\section{Problem setting and assumptions}\label{sec:defs}

To solve the optimization problem \eqref{eq:opt}, we consider iterative first-order methods described by dynamics of the form
\begin{align}\label{eq:ssalg}
		\xi^{t+1} &= A \xi^t + B (u^t + w^t), &
		y^t &= C \xi^t, &
		u^t &= \df( y^t ),
\end{align}
where $\xi^t \in \R^{n\times d}$ is the \textit{state} of the algorithm, $y^t\in\R^{1\times d}$ is where we evaluate the gradient, $u^t\in\R^{1\times d}$ is the (exact) gradient, $w^t\in\R^{1\times d}$ is the noise, and $t$ is the \textit{iteration}. The state is the \emph{memory} of the algorithm; its size reflects the number of past iterates that must be stored at each timestep. 
Given an algorithm $(A,B,C)$, objective function $f$, and initial condition $\xi^0$, a \emph{trajectory} is any sequence $(\xi^t, u^t, y^t, w^t)_{t\ge 0}$ that satisfies~\eqref{eq:ssalg}. This general framework encompasses a wide variety of fixed-parameter first-order iterative methods; we discuss this in more detail in \cref{sec:parameterization}.

\begin{remark}[important notational convention]
  Throughout this paper, we represent iterates as \emph{row vectors}, so the matrices $A\in\R^{n\times n}$, $B\in\R^{n\times 1}$, and $C\in\R^{1\times n}$ act on the \emph{columns} of the iterates while the gradient oracle $\df: \R^{1\times d} \to \R^{1\times d}$ acts on the \emph{rows}. Thus, $\norm{\cdot}$ denotes the Frobenius norm, for example, $\norm{\xi^t}^2 \defeq \sum_{i=1}^n\sum_{j=1}^d (\xi^t_{ij})^2$.
\end{remark}

A necessary condition for \eqref{eq:algform} to represent a valid iterative algorithm is that there exists a fixed point corresponding to the stationary point of the objective function. In other words, there should exist $\xi^\star \in \R^{n\times d}$ and $y^\star\in \R^{1\times d}$ satisfying
$\xi^\star = A\xi^\star$ and $y^\star = C\xi^\star$.
Taking this condition column-wise, we obtain the following.

\begin{assumption}[algorithm form]\label{ass:algform}
  The matrix $A$ has an eigenvalue at $1$, and the associated eigenvector $v\in\R^n$ satisfies $Cv\neq 0$.\footnote{In the language of control theory, the discrete-time system $(A,B,C)$ has an \emph{observable integrator}.}
\end{assumption}

\subsection{Performance evaluation}\label{sec:perf_eval}

Let $\alg = (A,B,C)$ denote an algorithm as in~\eqref{eq:ssalg} and let $\F \in \{Q_{m,L},F_{m,L}\}$ denote one of function classes defined in \cref{sec:intro}.

\parag{Convergence rate} The convergence rate describes the first phase of convergence observed in \cref{fig:demo_rand_trajectories}: exponential decrease of the error. In this regime, gradients are relatively large compared to the noise, so we set $w^t=0$ for all $t\geq 0$. For any algorithm $\alg$ with initial point $\xi^0$ and fixed point $\xi^\star$ and any function $f\in\F$, consider the trajectory $(\xi^0,\xi^1,\dots)$ produced by $\alg$. We define the linear convergence rate as
\begin{equation}\label{def:conv_rate}
	\rate(\alg,\F) \defeq \inf \left\{ r > 0 \;\Bigg\vert\; \sup_{f\in\F}\,\sup_{\xi^0\in\R^{n\times d}}\, \sup_{t\ge 0}\, \frac{ \norm{\xi^t-\xi^\star} }{ r^t \norm{\xi^0-\xi^\star} } < \infty \right\}.
\end{equation}
If the convergence rate $r = \rate(\alg,\F)$ is finite, then for all $\epsilon > 0$, there exists $c>0$ such that the trajectory satisfies $\norm{\xi^t-\xi^\star} \leq c (r+\epsilon)^t \norm{\xi^0-\xi^\star}$ for all functions $f\in\F$, initial points $\xi^0\in\R^{n\times d}$, and iterations $t\ge 0$.
This definition corresponds to the conventional notion of \emph{linear convergence rate} used in the worst-case analysis of deterministic algorithms. If $r < 1$, the algorithm is said to be \emph{globally linearly convergent}, and for all $\epsilon > 0$, we have $\norm{y^t-y^\star} = O((r+\epsilon)^t)$. Smaller $r$ corresponds to a faster (worst-case) convergence rate. 

\parag{Sensitivity} The sensitivity characterizes the steady-state phase of convergence observed in \cref{fig:demo_rand_trajectories}. The steady-state error depends on the noise characteristics. We make the following assumptions on the noise sequence.

\begin{assumption}[Noise sequence]\label{ass:noise}
  We assume that the noise sequence $w^0,w^1,\dots$ has joint distribution $\mathbb{P}\in\mathcal{P}_\sigma$, with parameter $\sigma$ to be defined shortly. We assume the set of admissible joint distributions $\mathcal{P}_\sigma$ satisfies:
  \begin{enumerate}
    \item \emph{Independence across time.} For all $\mathbb{P} \in \mathcal{P}_\sigma$, if $w\sim \mathbb{P}$, then $w^t$ and $w^\tau$ are independent for all $t\ne \tau$. Then we may characterize the joint distribution $\mathbb{P} \in \mathcal{P}_\sigma$ by its associated marginal distributions $(\mathbb{P}^0,\mathbb{P}^1,\dots)$. We do not assume the $\mathbb{P}^t$ are necessarily identical.
    \item \emph{Zero-mean and bounded variance.} For all $(\mathbb{P}^0,\mathbb{P}^1,\dots) \in \mathcal{P}_\sigma$, we have that $\E_{w^t\sim\mathbb{P}^t}(w^t) = 0$
    and $\E_{w^t\sim\mathbb{P}^t}(\norm{w^t}^2) \leq \sigma^2$. 
  \end{enumerate}
\end{assumption}
\Cref{ass:noise} implies $\mathcal{P}_\sigma$ is completely characterized by the variance bound~$\sigma^2$.
For algorithm $\alg$, function class $\F$, initial $\xi^0$, and family of distributions $\mathcal{P}_\sigma$, consider the stochastic iterate sequence $y^0,y^1,\dots$ produced by $\alg$.
We define noise sensitivity as\looseness=-1
\begin{equation}\label{def:gamma}
  \sensitivity(\alg,\F) \defeq \sup_{f\in\F}\, \sup_{\xi^0\in\R^{n\times d}}\, \sup_{\mathbb{P}\in \mathcal{P}_\sigma}\,\, \limsup_{T\to \infty}  \sqrt{\E_{w\sim \mathbb{P}}\, \frac{1}{\sigma^2T} \sum_{t=0}^{T-1} \norm{y^t-y^\star}^2}.
\end{equation}
For all cases we consider, the normalized definition of sensitivity \eqref{def:gamma} does not depend on $\sigma$.
A smaller sensitivity is desirable because it means the algorithm achieves small error in spite of gradient perturbations.
This definition is similar to that used in recent works exploring first-order algorithms
using techniques from robust control~\cite{mihailo,mert,scherer}.

\begin{remark}\label{rem:state_transformation}
  Some previous works have defined sensitivity with respect to the squared-norm of the state $\|\xi^t-\xi^\star\|^2$ \cite{scherer,mihailo}. This quantity, however, depends on the state-space realization; performing the coordinate transformation $(A,B,C) \mapsto (TAT^{-1}, TB, CT^{-1})$ for some invertible matrix $T$ does not change the sequence of inputs $u^t$ or outputs $y^t$, but it transforms the states: $\xi^t\mapsto T \xi^t$. Our definition of sensitivity \eqref{def:gamma} is invariant under such transformations.
\end{remark}

\subsection{Algorithm parameterization}\label{sec:parameterization}

While our \textit{analysis} results apply to the general model~\eqref{eq:ssalg}, for  \emph{design} we will further restrict the class of algorithms to those with state dimension $n=2$. Algorithms of the form~\eqref{eq:ssalg} have $n^2+2n$ degrees of freedom in the matrices $A,B,C$, so the case $n=2$ should require $8$ parameters. However, many of these parameters are redundant, and under \cref{ass:algform}, it turns out the case $n=2$ is completely characterized by the three-parameter family~\eqref{eq:algform0}.

\begin{lemma}[three-parameter family]\label{lem:alg3}
  Under \cref{ass:algform}, 
  any algorithm of the form~\eqref{eq:ssalg} with state dimension $n=2$ 
  is equivalent to an algorithm in the three-parameter family~\eqref{eq:algform0} for some $(\alpha,\beta,\eta)$. By equivalent, we mean that both algorithms produce the same sequence of iterates $(y^0,y^1,\ldots)$ given appropriate initialization.
\end{lemma}
\begin{proof}
    See \cref{app:proof_3_param}
\end{proof}

We refer to specific algorithms from the three-parameter family~\eqref{eq:algform0} by their parameters $(\alpha,\beta,\eta)$. All such algorithms satisfy \cref{ass:algform}. For specific values for the parameters, the algorithm~\eqref{eq:algform0} recovers several known algorithms, shown in \cref{table:algos}, which will serve as useful benchmarks for our designs. Later in \cref{sec:justification_three_param}, we provide numerical evidence that further justifies our choice of parameterization.

\begin{table}[ht]
	\centering
	\caption{Comparison of different algorithms with their recommended/standard tunings. For RM, the parameter satisfies $1-\sqrt{\frac{m}{L}} \leq r \leq 1-\frac{m}{L}$, and interpolates between TM and GD with $\alpha=\frac{1}{L}$.} 
  \label{table:algos}
	\begin{tabular}[t]{lcccc}
		\toprule
		Algorithm name & $\alpha$ & $\beta$ & $\eta$ \\
		\midrule
		Gradient Descent (GD) & $\frac{1}{L}$ or $\frac{2}{L+m}$ & 0 & 0 \\
		Heavy Ball (HB), \cite{polyak_book} & $\frac{4}{(\sqrt{L}+\sqrt{m})^2}$ & $\left(\frac{\sqrt{L}-\sqrt{m}}{\sqrt{L}+\sqrt{m}}\right)^2$ & 0 \\
		Fast Gradient (FG), \cite{nesterov_book} & $\frac{1}{L}$ & $\frac{\sqrt{L}-\sqrt{m}}{\sqrt{L}+\sqrt{m}}$ & $\frac{\sqrt{L}-\sqrt{m}}{\sqrt{L}+\sqrt{m}}$ \\
		Triple Momentum (TM), \cite{tmm} & $\frac{2\sqrt{L}-\sqrt{m}}{L^{3/2}}$ & $\frac{(\sqrt{L}-\sqrt{m})^2}{L+\sqrt{mL}}$ & $\frac{(\sqrt{L}-\sqrt{m})^2}{2L-m+\sqrt{mL}}$ \\
		Robust Momentum (RM), \cite{rmm} & $\frac{(1-r)^2(1+r)}{m}$ & $\frac{Lr^3}{L-m}$ & $\frac{mr^3}{(L-m)(1-r)^2(1+r)}$ \\
		\bottomrule
	\end{tabular}
\end{table}

Different triples $(\alpha,\beta,\eta)$ generally correspond to different algorithms, with one important exception: Gradient Descent has a degenerate family of parameterizations.

\begin{proposition}\label{prop:GD_unique_param}
	Gradient Descent with parameterization $(\alpha,0,0)$ can also be parameterized by $\bigl(\alpha(1-\beta),\,\beta,\,\frac{\beta}{1-\beta}\bigr)$ for any choice of $\beta \neq 1$.
\end{proposition}

To see why this holds, substitute $(\alpha,\beta,\eta)\mapsto\bigl(\alpha(1-\beta),\beta,\frac{\beta}{1-\beta}\bigr)$ into~\eqref{eq:algform0}, and observe that the update becomes GD applied to the quantity $y^t \defeq \frac{1}{1-\beta}(x^t - \beta x^{t-1})$.

\begin{remark}
    \cref{lem:alg3,prop:GD_unique_param} can also be proved using the notion of a \emph{transfer function} (TF) \cite[\S2.7]{antsaklis}, which is the linear map from the $z$-transform of $(u^t + w^t)$ to the $z$-transform of $y^t$. The TF of \eqref{eq:ssalg} is $G(z) = C(zI-A)^{-1}B$ and when two dynamical systems have the same TF, they are equivalent in the sense of \cref{lem:alg3}. When $n=2$, $G(z)$ is a rational function with one zero and two poles, one of which must be at $z=1$ by \cref{ass:algform}. This leaves three degrees of freedom (the second pole, the zero, and a constant gain). Similarly, it is easy to check that the TF for $\bigl(\alpha(1-\beta),\,\beta,\,\frac{\beta}{1-\beta}\bigr)$ is $G(z) = \frac{-\alpha}{z-1}$, which is independent of $\beta$.
\end{remark}

In the next two sections, we focus on $Q_{m,L}$ and $F_{m,L}$. For each, we provide a tractable approach for bounding the convergence rate and sensitivity, and we design algorithms of the form~\eqref{eq:algform0} that effectively trade off these performance metrics.

\section{Strongly convex quadratic functions}\label{sec:Q}

\subsection[Performance bounds]{Performance bounds for \texorpdfstring{$Q_{m,L}$}{Q}}

Quadratics have been treated extensively in recent works. For this class, the convergence rate has been characterized~\cite{lessard16}, and closed-form expressions for the sensitivity of GD, HB, and FG have been obtained~\cite{mihailo,mert}. We now present versions of these results adapted to our algorithm class of interest.

\begin{lemma}[$Q_{m,L}$ analysis]\label{prop:perf_quad_general}
Consider algorithm $\alg = (A,B,C)$ defined in~\eqref{eq:ssalg} satisfying \cref{ass:algform} applied to $f \in Q_{m,L}$ defined in \cref{sec:intro}. Assume the noise sequence satisfies \cref{ass:noise}. Algorithm $\alg$ has convergence rate
	\[
		\rate(\alg,Q_{m,L}) = \sup_{q \in [m,L]}\rho\bigl( A + q B C \bigr),
	\]
    where $\rho(\cdot)$ denotes the spectral radius of a matrix (largest eigenvalue magnitude).
	If $\rate(\alg,Q_{m,L}) < 1$, the algorithm has finite sensitivity given by
	\[
		\sensitivity(\alg,Q_{m,L})
		= \sup_{q\in[m,L]} \sqrt{B^\tp P_q B},
	\]
	where $P_{q}$ is the unique solution to the linear equation
    \begin{equation}\label{eq:lyap-eqn}
        (A+q BC)^\tp P_{q} (A+q BC) - P_{q} + C^\tp C = 0.
    \end{equation}
\end{lemma}

\begin{proof}
    See \cref{app:proof_QmL}.
\end{proof}

The linear equation~\eqref{eq:lyap-eqn} is a \textit{Lyapunov equation}, and the fact that it has a unique solution is provided in \cite[\S11.4.4]{fadali2013digital}.
Similar results have appeared in the context of algorithm analysis for quadratics~\cite{mert,mihailo} and use that the sensitivity is equivalent to the $\mathcal{H}_2$ norm of the system, which can be computed using a Lyapunov approach.

The expressions for the rate and sensitivity in \cref{prop:perf_quad_general} may be difficult to evaluate if the matrices $(A,B,C)$ are large\footnote{Neither $\rho(A+q B C)$ nor $P_q$ are convex functions of $q$ in general.}. Fortunately, the sizes of these matrices only depend on the state dimension $n$, which is typically small ($n\leq 2$ for all methods in \cref{table:algos}).
Moreover, the dimension $d$ of the function domain does not affect the complexity of the convergence rate or sensitivity formula.

\subsection[Algorithm design]{Algorithm design for $Q_{m,L}$}\label{sec:algodesign_QmL}

For strongly convex quadratic functions, first-order methods can achieve exact convergence in $d$ iterations when there is no gradient noise, where $d$ is the dimension of the domain of $f$. One such example is the Conjugate Gradient (CG) method~\cite[Thm.~5.4]{nocedal_wright}.
However, when the number of iterations $t$ satisfies $t < d$, exact convergence is not possible in general. Nesterov's lower bound~\cite[Thm.~2.1.13]{nesterov_book} demonstrates that for any $t \ge 0$, one can construct a function $f\in Q_{m,L}$ with domain dimension $d > t$ such that:
\begin{equation}\label{eq:nesterov_lowerbound}
	\norm{y^t-y^\star} \ge \biggl( \frac{ \sqrt{L}-\sqrt{m}}{\sqrt{L}+\sqrt{m}} \biggr)^t \norm{y^0-y^\star}.
\end{equation}
This lower bound holds for any first-order method such that $y^t$ is a linear combination of $y^0$ and (the exact) past gradients $\df(y^0),\dots,\df(y^{t-1})$. This class includes not only CG but also methods with unbounded memory.

In the regime $t < d$, the CG method matches Nesterov's lower bound~\cite[Thm.~5.5]{nocedal_wright} and is therefore optimal in terms of worst-case rate. However, it is not clear how CG should be adjusted to be robust in the presence of additive gradient noise, since it has no tunable parameters.

The HB method~\eqref{eq:hb}, when tuned as in \cref{table:algos}, also matches Nesterov's lower bound when applied to the function class $Q_{m,L}$~\cite[\S3.2.1]{polyak_book}, but has a simpler implementation than CG: its updates are linear and its parameters are constant.

We adopted the three-parameter class $(\alpha,\beta,\eta)$ described in \cref{sec:parameterization} as our search space for optimized algorithms because HB belongs to this class and achieves optimal performance on $Q_{m,L}$ when there is no noise. Substituting the three-parameter algorithm~\eqref{eq:algform} into \cref{prop:perf_quad_general}, we obtain the following result.

\begin{corollary}[$Q_{m,L}$ analysis, reduced]\label{cor:perf_quad}
	Consider the three-parameter algorithm $\alg=(\alpha,\beta,\eta)$ defined in \eqref{eq:algform0} applied to a strongly convex quadratic $f\in Q_{m,L}$ defined in \cref{sec:intro}, and assume the noise sequence satisfies \cref{ass:noise}. The algorithm has convergence rate
  \begin{multline}\label{eq:rho_mL_quad}
		\rate(\alg,Q_{m,L})
		= \max_{q \in \{m,L\}}
    \begin{cases}
    	\sqrt{\beta - \alpha \eta q} & \text{if }\Delta < 0, \\
    	\frac{1}{2} \left( \lvert \beta+1-\alpha q-\alpha\eta q \rvert + \sqrt{\Delta}\right) & \text{if }\Delta \geq 0,
    \end{cases}\\
    \text{where }\Delta\defeq (\beta+1-\alpha q-\alpha\eta q)^2-4(\beta-\alpha\eta q).
	\end{multline}
	If $\rate(\alg,Q_{m,L}) < 1$, the algorithm has sensitivity
    	\begin{equation}\label{eq:gamma_mL_quad} 
	\sensitivity(\alg,Q_{m,L})
		= \max_{q\in\{m,L\}} \sqrt{h(q)} \defeq \sqrt{\frac{\alpha(1+\beta+(1+2\eta)\alpha\eta q) }{q(1-\beta+\alpha\eta q)(2+2\beta -(1+2\eta)\alpha q)}}
	\end{equation}
\end{corollary}

\begin{proof}
See \cref{app:proof_QmL_reduced}.
\end{proof}

In \cref{cor:perf_quad}, the suprema from \cref{prop:perf_quad_general} are replaced by a simple maximum; we only need to check the endpoints of the interval $[m,L]$.

We begin with a baseline; the convergence and sensitivity of GD on $Q_{m,L}$.

\begin{theorem}[GD analysis for $Q_{m,L}$]\label{thm:RGM}
	Consider the function class $Q_{m,L}$ and let $r$ be a parameter chosen with $\frac{L-m}{L+m} \le r < 1$. Then, under \cref{ass:noise}, the algorithm $\alg$ of the form~\eqref{eq:algform0} with
	$\alpha = \frac{1}{m}(1-r)$, 
	$\beta = 0$, and $\eta = 0$
	achieves
	\[
		\rate(\alg, Q_{m,L}) = r,\qquad\text{and}\qquad
		\sensitivity(\alg,Q_{m,L}) 
        = \frac{1}{m}\sqrt{\frac{1-r}{1+r}}.
	\]
\end{theorem}
\begin{proof}
    Omitted; follows from straightforward algebra.
\end{proof}

Our proposed algorithm for the class $Q_{m,L}$ is a special tuning of Heavy Ball, which we named \emph{Robust Heavy Ball} (RHB). The RHB algorithm was found by careful analysis of the analytic expressions for the rate and sensitivity in \cref{cor:perf_quad}.
The algorithm is described in the following theorem.

\begin{theorem}[Robust Heavy Ball, RHB]\label{thm:RHB}
	Consider the function class $Q_{m,L}$ and let $r$ be a parameter chosen with $\tfrac{\sqrt{L}-\sqrt{m}}{\sqrt{L}+\sqrt{m}} \leq r < 1$. Then under \cref{ass:noise}, the algorithm $\alg$ of the form~\eqref{eq:algform0} with
	$\alpha = \frac{1}{m}(1-r)^2$,
	$\beta = r^2$, and $\eta=0$
	achieves
	\[
		\rate(\alg,Q_{m,L}) = r \qquad\text{and}\qquad
		\sensitivity(\alg,Q_{m,L})
		= \frac{1}{m}\sqrt{\frac{1-r^4}{(1+r)^4}}.
	\]
\end{theorem}

\begin{proof}
  We first show that the parameter $r$ is in fact the convergence rate of the algorithm. Substituting the algorithm parameters $\alpha = \frac{1}{m}(1-r)^2$, $\beta = r^2$, and $\eta = 0$ into~\eqref{eq:rho_mL_quad}, we obtain
  $
    \Delta = -(1-r)^4 \left(\frac{q }{m}-1\right)
    \bigl(\bigl(\frac{1+r}{1-r}\bigr)^2-\frac{q}{m}\bigr)
  $.
  Rearranging the inequalities $m\leq q \leq L$ and $\tfrac{\sqrt{L}-\sqrt{m}}{\sqrt{L}+\sqrt{m}} \leq r < 1$ yields $1 \leq \frac{q}{m} \leq \frac{L}{m} \leq \bigl(\frac{1+r}{1-r}\bigr)^2 < \infty$.
  Thus, we conclude that $\Delta \leq 0$ and we have $\rate(\alg,Q_{m,L}) = \sqrt{\beta} = r$, as required.
  Now substitute the parameters into \eqref{eq:gamma_mL_quad} and obtain
  \[
    h(m)-h(L)
    =
    \frac{(1-r) (r^2+1) (\frac{L}{m}-1) \bigl(\bigl(\frac{1+r}{1-r}\bigr)^2-\frac{L}{m}\bigr)}{Lm (r+1)^3 \bigl(1+\bigl(\frac{1+r}{1-r}\bigr)^2-\frac{L}{m}\bigr)} \geq 0,
  \]
  so $h(q)$ is maximized when $q=m$, and $\sensitivity(\alg,Q_{m,L}) = \sqrt{h(m)} = \frac{1}{m}\sqrt{\frac{1-r^4}{(1+r)^4}}$.
\end{proof}

Although we set out to design an algorithm in the three-parameter class $(\alpha,\beta,\eta)$, our designed algorithm RHB uses ${\eta=0}$, so it is a particular tuning of HB. Our numerical experiments in \cref{sec:pareto_validation} suggest that this parameter choice yields the most effective trade-off between rate and sensitivity. In other words, it is unnecessary to use a nonzero $\eta$ when optimizing over the class $Q_{m,L}$.

If we choose the smallest possible $r = \frac{\sqrt{L}-\sqrt{m}}{\sqrt{L}+\sqrt{m}}$ (fastest possible convergence rate), then we recover Polyak's tuning of HB, whose convergence rate matches Nesterov's lower bound. It is straightforward to check that the sensitivity is a monotonically decreasing function of $r$, so as the convergence rate slows down ($r$ increases), the algorithm becomes less sensitive to noise. Thus, $r$ is a single tunable parameter that enables us to trade off convergence rate with sensitivity to noise.

\section{Smooth strongly convex functions}\label{sec:F}

We now consider the class $F_{m,L}$ of strongly convex functions whose gradient is Lipschitz continuous.
In contrast to the $Q_{m,L}$ case, the function class $F_{m,L}$ is not readily parameterizable. Therefore, we use a Lyapunov approach to certify performance bounds.

\subsection{Lyapunov analysis}\label{sec:lyapunov}

Our analysis is based on searching for a function that certifies either a particular convergence rate (assuming no noise) or level of sensitivity (assuming noise). When used in the context of certifying stability of dynamical systems, such a function is called a \textit{Lyapunov} function. This function will depend on the state $\x^t\in X$ of a (to-be-defined) dynamical system. In \cref{sec:lift}, we will show how to construct this augmented system from the algorithm~\eqref{eq:ssalg}. Before doing so, we first show how to use such a function to certify performance.

\paragraph{Bounds on the rate of convergence}

To bound the rate of convergence of an algorithm, we search for a function $V(\x)$ that decreases along trajectories and is lower bounded by the squared norm of the iterates.

\begin{lemma}[Lyapunov analysis for rate of convergence]\label{lem:lyap_rate}
  Consider an algorithm $\alg = (A,B,C)$ defined in \eqref{eq:ssalg} satisfying \cref{ass:algform} applied to a function $f\in\F$ with no gradient noise ($w^t=0$ for all $t$). If there exists a function $V : X\to\R$ and a constant $r>0$ such that, for all functions $f\in\F$ and all iterations $t\geq 0$, the iterates $\x^t$ satisfy the conditions
  \begin{enumerate}[(i)]
    \item \emph{Lower bound condition:} $V(\x^t) \geq \|\xi^t-\xi^\star\|^2$ and
    \item \emph{Decrease condition:} $V(\x^{t+1}) \leq r^2 \,V(\x^t)$,
  \end{enumerate}
  then $\rate(\alg,\F)\leq r$.
\end{lemma}

\begin{proof}
  By applying the lower bound condition followed by the decrease condition at each iteration $t\geq 0$, we obtain the chain of inequalities
$\|\xi^t-\xi^\star\|^2 \leq V(\x^t) \leq r^2\,V(\x^{t-1}) \leq \ldots \leq r^{2t}\,V(\x^0)$.
  From \cref{def:conv_rate}, we obtain $\rate(\alg,\F)\leq r$, as required.
\end{proof}

\paragraph{Bounds on the sensitivity to noise}

To bound the sensitivity of an algorithm to noise, we search for a function $V(\x)$ that satisfies the following conditions.

\begin{lemma}[Lyapunov analysis for sensitivity to noise]\label{lem:lyap_sensitivity}
  Consider an algorithm $\alg = (A,B,C)$ defined in \eqref{eq:ssalg} satisfying \cref{ass:algform} applied to a function $f\in\F$ with additive gradient noise satisfying \cref{ass:noise}. If there exists a function $V : X\to\R$ and a constant $\gamma>0$ such that, for all functions $f\in\F$ and all iterations $t\geq 0$, the iterates $\x^t$ satisfy the conditions
  \begin{enumerate}[(i)]
    \item \emph{Lower bound condition:} $\E V(\x^t) \geq 0$ and
    \item \emph{Decrease condition:} $\E V(\x^{t+1}) - \E V(\x^t) + \E \norm{y^t-y^\star}^2 \leq \sigma^2\gamma^2$,
  \end{enumerate}
  then $\sensitivity(\alg,\F)\leq\gamma$.
\end{lemma}

\begin{proof}
  Averaging the decrease condition over $t=0,\dots,T-1$, we obtain
  \[
    \frac{1}{T} \E V(\x^T) - \frac{1}{T} \E V(\x^0) + \frac{1}{T} \sum_{t=0}^{T-1} \E \norm{y^t-y^\star}^2 \leq \sigma^2\gamma^2.
  \]
  Applying the lower bound condition to the first term and taking the limit superior as $T\to\infty$, the second term vanishes and the result follows from \eqref{def:gamma}.
\end{proof}

In both cases (\cref{lem:lyap_rate,lem:lyap_sensitivity}), we will see that searching over functions $V(\x)$ of a particular form can be cast as a semidefinite program.

\subsection{Lifted dynamics}\label{sec:lift}

To obtain performance bounds, we use a \emph{lifting} approach that searches for certificates of performance that depend on a finite history of past algorithm iterates and function values. The main idea is to lift the state to a higher dimension so that we can search over a broader class of certificates to reduce the conservativeness of the bound. We denote the \textit{lifting dimension} by $\ell\ge 0$, which dictates the dimension of the lifted state, and we use boldface symbols to denote quantities related to the lifted dynamics.

Given a trajectory of the system~\eqref{eq:ssalg}, we define the following augmented vectors, each consisting of $\ell+1$ consecutive iterates of the system:
\begin{equation}\label{eq:lifted_vars}
  \y^t \defeq \bmat{ y^t-y^\star \\ \vdots \\  y^{t-\ell}-y^\star }, \qquad
  \u^t \defeq \bmat{ u^t \\ \vdots \\  u^{t-\ell} }, \qquad\text{and}\qquad
  \f^t \defeq \bmat{ f^t-f^\star \\ \vdots \\  f^{t-\ell}-f^\star },
\end{equation}
where $\y^t,\u^t\in\R^{(\ell+1)\times d}$ and $\f^t\in\R^{\ell+1}$ with $f^t \defeq f(y^t)$ and $f^\star \defeq f(y^\star)$ (recall our convention that algorithm inputs and outputs $u^t$, $y^t$ are \emph{row} vectors). Also, define the truncation matrices $Z,Z_+ \in \R^{\ell \times (\ell+1)}$ as
\begin{equation}\label{eq:Z}
  Z_+ \defeq \bmat{I_\ell & 0_{\ell\times 1}}
  \quad\text{and}\quad
  Z \defeq \bmat{0_{\ell\times 1} & I_\ell}.
\end{equation}
Multiplying an augmented vector on the left by $Z$ removes the most recent iterate at time $t$, while multiplication by $Z_+$ removes the last iterate at time $t-\ell$. Using these augmented vectors, we then define the augmented state
\begin{equation}\label{eq:state_aug}
	\x^t \defeq (\bxi^t,Z \f^t), \qquad\text{where}\qquad
    \bxi^t \defeq \bmat{ \xi^t-\xi^\star \\ Z \y^t \\ Z \u^t} \in \R^{(n+2\ell)\times d},
\end{equation}
which consists of the deviations of the state $\xi^t$ and past inputs $y^{t-1},\ldots,y^{t-\ell}$, outputs $u^{t-1},\ldots,u^{t-\ell}$, and function values $f^{t-1},\ldots,f^{t-\ell}$ of the original system from equilibrium.
The augmented dynamics, which follow from \cref{eq:ssalg,eq:lifted_vars,eq:Z,eq:state_aug}, are
\begin{subequations}\label{eq:ss_aug}
  \begin{align}
    \bxi^{t+1} &= \underbrace{\bmat{ A & 0 & 0 \\ Z_+ \e_1 C & Z_+ Z^\tp & 0 \\ 0 & 0 & Z_+ Z^\tp }}_{\A} \bxi^t + \underbrace{\bmat{B \\ 0 \\ Z_+ \e_1}}_{\B} u^t + \underbrace{\bmat{B \\ 0 \\ 0}}_{\H} w^t, \\
    \bmat{\y^t \\ \u^t} &= \underbrace{\bmat{ \e_1 C & Z^\tp & 0 \\ 0 & 0 & Z^\tp }}_{\C} \bxi^t + \underbrace{\bmat{ 0 \\ \e_1}}_{\D} u^t,
  \end{align}
\end{subequations}
where $\e_1=(1,0,\ldots,0)\in\R^{\ell+1}$. We can recover the iterates of the original system by projecting the augmented state and the input as
\begin{align}
	\xi^t-\xi^\star \!=\! \underbrace{\squeezemat\bmat{I_n & 0_{n\times(2\ell+1)}}}_{\X}\! \bmat{\bxi^t \\ u^t}\!, \
	y^t-y^\star \!=\! \underbrace{\squeezemat\bmat{C & 0_{1\times(2\ell+1)}}}_{\Y}\! \bmat{\bxi^t \\ u^t}\!, \
	u^t \!=\! \underbrace{\squeezemat\bmat{0_{1\times(n+2\ell)} & 1}}_{\U}\! \bmat{\bxi^t \\ u^t}\!. \label{eq:ss_aug_proj}
\end{align}

\begin{remark}[State reduction for noise-free case]
When there is no noise (${w^t=0}$), the component $\bxi^t$ of the augmented state~\eqref{eq:state_aug} has linearly dependent rows; knowledge of the past state $\xi^{t-\ell}$ and subsequent inputs $u^{t-\ell},\dots,u^{t-1}$ is enough to reconstruct the outputs $y^{t-\ell},\dots,y^{t-1}$. Thus, in this case, we could work with a smaller lifted state that does not include outputs. Such a reduction does not change the results of the analysis, but makes the semidefinite programs smaller and therefore (slightly) more computationally efficient.
\end{remark}

\subsection{Interpolation conditions}

A useful characterization of the function class $F_{m,L}$ is given by the \emph{interpolation conditions}. These are necessary and sufficient conditions on a sequence of points to be interpolable by a function in the class.
We state the result from~\cite[Thm.~4]{taylor2017smooth}, rephrased to match our notation.

\begin{proposition}[Interpolation conditions for $F_{m,L}$]
	\label{prop:interp}
	Let $y^1,\dots,y^k \in \R^{1\times d}$ and $u^1,\dots,u^k \in \R^{1\times d}$ and $f^1,\dots, f^k \in \R$. The following statements are equivalent.

	\begin{enumerate}[label=\arabic*)]
		\item There exists $f\in F_{m,L}$ such that $f(y^i) = f^i$ and $\df(y^i) = u^i$ for $i=1,\dots,k$.
		\item For all $i,j \in \{1,\dots,k\}$, the following inequality holds:
		\begin{equation}\label{eq:interp}
			\trace\bigl((u^i)^\tp (y^i-y^j)\bigr) - (f^i-f^j) 
      + \frac{1}{2\,(L-m)}
			\trace\bmat{ y^i-y^j \\ u^i-u^j }^\tp
			\bmat{ -mL & m \\ m & -1 }
			\bmat{ y^i-y^j \\ u^i-u^j } \geq 0.
		\end{equation}
	\end{enumerate}
\end{proposition}

We now develop a version of \cref{prop:interp} with a single parameterized inequality that holds for the augmented vectors defined in~\eqref{eq:lifted_vars}.

\begin{lemma}\label{lem:cvx_func_ineq}
  Consider an algorithm \eqref{eq:ssalg} applied to a function $f \in F_{m,L}$. Define the augmented iterates in \eqref{eq:lifted_vars} and the index set $I \defeq \{1,\ldots,\ell+1,\star\}$, and let $\e_i$ denote the $i\textsuperscript{th}$ unit vector in $\R^{\ell+1}$ with $\e_\star \defeq 0\in\R^{\ell+1}$. Then the inequality
  \begin{equation}\label{eq:cvx_ineq}
    \trace \bmat{\y^t \\ \u^t}^\tp \Pi(\Lambda) \bmat{\y^t \\ \u^t} + \pi(\Lambda)^\tp \f^t \ge 0
  \end{equation}
  holds for all $\Lambda\in\R^{(\ell+2)\times(\ell+2)}$ such that $\Lambda\geq 0$ (elementwise), where
  \begin{subequations}\label{eq:multipliers}
  \begin{align}
    \Pi(\Lambda) &\defeq \sum_{i,j\in I} \Lambda_{ij} \bmat{ -mL\,(\e_i-\e_j)(\e_i-\e_j)^\tp & (\e_i-\e_j)(m\e_i-L\e_j)^\tp \\[2pt] (m\e_i-L\e_j)(\e_i-\e_j)^\tp & -(\e_i-\e_j)(\e_i-\e_j)^\tp }, \\
    \pi(\Lambda) &\defeq 2\,(L-m) \sum_{i,j\in I} \Lambda_{ij}\,(\e_i-\e_j).
  \end{align}
  \end{subequations}
\end{lemma}

\begin{proof}
  For each $i\in I$, define the vectors $\tilde y_i \defeq \e_i^\tp \y^t$ and $\tilde u_i \defeq \e_i^\tp \u^t$ and $\tilde f_i \defeq \e_i^\tp \f^t$. By definition, the points $(\tilde y_i,\tilde u_i,\tilde f_i)$ are interpolated by the function $f\in F_{m,L}$, so by \cref{prop:interp} the interpolation conditions \eqref{eq:interp} are satisfied. The proof is then completed by noting that the inequality \eqref{eq:cvx_ineq} can be expanded as
  \[
    \sum_{i,j\in I} \Lambda_{ij} \Bigl(-mL \|\tilde y_i-\tilde y_j\|^2 + 2\,(\tilde y_i-\tilde y_j) (m \tilde u_i-L \tilde u_j)^\tp - \|\tilde u_i - \tilde u_j\|^2 + 2 (L-m) (\tilde f_i-\tilde f_j)\Bigr) \geq 0,
  \]
  which is a nonegative weighted combination of the interpolation conditions, where the interpolation condition between index $i$ and $j$ is scaled by $2\,(L-m)\Lambda_{ij}\geq 0$.
\end{proof}

\subsection[Performance bounds]{Performance bounds for \texorpdfstring{$F_{m,L}$}{F}}\label{sec:perf_cvx}

We now use the interpolation conditions to search for a Lyapunov function that depends on the state $\x^t$ of the lifted system that satisfies the conditions in \cref{sec:lyapunov} to certify either the convergence rate or sensitivity.
Motivated by the fact that the nonnegative inequalities in~\eqref{eq:cvx_ineq} are quadratic in the inputs and outputs of the gradient and linear in the function values, we search for certificates $V : \R^{(n+2\ell)\times d}\times\R^\ell\to\R$ of the form
\begin{equation}\label{eq:lyap}
  V(\x) \defeq \trace\bigl(\bxi^\tp P \bxi\bigr) + p^\tp \f,
\end{equation}
where $\x = (\bxi,\f)$ is the state of the lifted system, and
we search over parameters $P=P^\tp \in\R^{(n+2\ell)\times(n+2\ell)}$ and  $p\in\R^\ell$.
Since $V$ is quadratic in the augmented state and linear in the augmented function values, we can efficiently search for such Lyapunov functions using the following linear matrix inequalities that leverage the characterization of smooth strongly convex functions in \cref{lem:cvx_func_ineq}.

\begin{theorem}[$F_{m,L}$ analysis]\label{thm:lmi_cvx}
	Consider an algorithm $\alg = (A,B,C)$ defined in \eqref{eq:ssalg} satisfying \cref{ass:algform} applied to a function $f \in F_{m,L}$ defined in \cref{sec:intro} with additive gradient noise satisfying \cref{ass:noise}.
    Define the matrices in \cref{eq:Z,eq:ss_aug,eq:ss_aug_proj,eq:multipliers}.

	If there exist $r > 0$, $P=P^\tp\in\R^{(n+2\ell)\times (n+2\ell)}$, $p\in\R^\ell$, and $\Lambda_1,\Lambda_2\geq 0$ such that
 	\begin{subequations}\label{eq:lmi_cvx_rate}
	\begin{align}
    \bmat{ \A & \B \\ I & 0 }^\tp \bmat{ P & 0 \\ 0 & -r^2 P } \bmat{ \A & \B \\ I & 0 } + \bmat{ \C & \D }^\tp \Pi(\Lambda_1) \bmat{ \C & \D } &\preceq 0
    \label{eq:lmi_cvx_1} \\
    (Z_+ - r^2 Z)^\tp p + \pi(\Lambda_1) &\leq 0
    \label{eq:lmi_cvx_2} \\
    \X^\tp \X - \bmat{I & 0}^\tp P \bmat{I & 0} + \bmat{ \C & \D }^\tp \Pi(\Lambda_2) \bmat{ \C & \D } &\preceq 0
    \label{eq:lmi_cvx_3} \\
    -Z^\tp p + \pi(\Lambda_2) &\leq 0
    \label{eq:lmi_cvx_4}
	\end{align}
    \end{subequations}
	then $\rate(\alg,F_{m,L}) \leq r$.
    
	If there exist $P=P^\tp\in\R^{(n+2\ell)\times(n+2\ell)}$, $p\in\R^\ell$, and $\Lambda_1,\Lambda_2\geq 0$ such that
	\begin{subequations}\label{eq:lmi_cvx_sensitivity}
	\begin{align}
    \bmat{ \A & \B \\ I & 0 }^\tp \bmat{ P & 0 \\ 0 & -P } \bmat{ \A & \B \\ I & 0 } + \bmat{ \C & \D }^\tp \Pi(\Lambda_1) \bmat{ \C & \D } + \Y^\tp \Y &\preceq 0
    \label{eq:lmi_cvx_sensitivity_1} \\
    (Z_+ - Z)^\tp p + \pi(\Lambda_1) &\leq 0
    \label{eq:lmi_cvx_sensitivity_2} \\
    -\bmat{I & 0}^\tp P \bmat{I & 0} + \bmat{ \C & \D }^\tp \Pi(\Lambda_2) \bmat{ \C & \D } &\preceq 0
    \label{eq:lmi_cvx_sensitivity_3} \\
    -Z^\tp p + \pi(\Lambda_2) &\leq 0
    \label{eq:lmi_cvx_sensitivity_4}
  \end{align}
  \end{subequations}
	then $\sensitivity(\alg,F_{m,L}) \leq \sqrt{\H^\tp P \H}$.
\end{theorem}

\begin{proof}
    See \cref{app:proof_FmL_analysis}.
\end{proof}

\begin{remark}\label{rem:ellzero}
  When the lifting dimension $\ell$ is zero, the lifted system is identical to the original system. In other words, $\xi^t = \bxi^t$. The system matrices also satisfy $A = \A$, and similarly for $B$. 
  As $\ell$ is increased, the LMIs~\eqref{eq:lmi_cvx_rate} and \eqref{eq:lmi_cvx_sensitivity} have the potential to yield less conservative bounds on the rate and sensitivity, respectively.
\end{remark}

Both bounds for the class $F_{m,L}$ in \cref{thm:lmi_cvx} can be evaluated and optimized efficiently. The sizes of the LMIs depend only on $n$ and $\ell$, which are typically small.

\subsection[Algorithm design]{Algorithm design for $F_{m,L}$}\label{sec:algodesign_FmL}

For the case with no noise, the Triple Momentum (TM) method~\cite{tmm} attains the fastest-known worst-case convergence rate of $1 - \sqrt{\frac{m}{L}}$ over the function class $F_{m,L}$. Recent work by Drori and Taylor \cite{oracle-complexity,ITEM} has confirmed that this rate is in fact optimal.

We adopted the three-parameter class $(\alpha,\beta,\eta)$ described in \cref{sec:parameterization} as our search space for optimized algorithms, because it includes TM as a special case, and FG, which is a popular choice for this function class.
We begin with the GD baseline.

\begin{theorem}[GD analysis for $F_{m,L}$]\label{thm:RGM2}
	Consider the function class $F_{m,L}$ and let $r$ be a parameter chosen with $\frac{L-m}{L+m} \le r < 1$. Then, under \cref{ass:noise}, the algorithm $\alg$ of the form~\eqref{eq:algform0} with
	$\alpha = \frac{1}{m}(1-r)$, 
	$\beta = 0$, and $\eta = 0$
	achieves
	\[
		\rate(\alg, Q_{m,L}) = r,\qquad\text{and}\qquad
		\sensitivity(\alg,Q_{m,L}) 
        = \frac{1}{m}\sqrt{\frac{1-r}{1+r}}.
	\]
\end{theorem}
\begin{proof}
See \cref{app:proof_GD_FmL}.
\end{proof}

GD achieves the same worst-case rate and sensitivity in $F_{m,L}$ (\cref{thm:RGM2}) as in $Q_{m,L}$ (\cref{thm:RGM}). This confirms the common knowledge that quadratics are worst-case function for GD.

Our proposed algorithm, the \emph{Robust Accelerated Method} (RAM), uses a parameter $r$ to trade off convergence rate and sensitivity to noise.
To design RAM, we applied the procedure outlined in~\cref{sec:design} to the rate LMI~\eqref{eq:lmi_cvx_rate} with lifting dimension $\ell=1$.
Our main result is \cref{thm:RAM}, which provides the exact convergence rate of RAM. While we do not construct a bound on the sensitivity, we show in \cref{sec:numerical} that it effectively trades off rate and sensitivity. To obtain a numerical bound on the sensitivity for a given parameter $r$, one can use the analysis result in \cref{thm:lmi_cvx}.

\begin{theorem}[Robust Accelerated Method, RAM]\label{thm:RAM}
	Consider the function class $F_{m,L}$, and let $r$ be a parameter chosen with $1-\sqrt{\frac{m}{L}} \leq r < 1$. Then, the algorithm $\alg$ of the form~\eqref{eq:algform0} with tuning
    {\[
    \alpha = \frac{(1+r)(1-r)^2}{m}, \quad
    \beta  = r\,\frac{L\,(1-r+2r^2)-m\,(1+r)}{(L-m)(3-r)},\quad
	\eta   = r\,\frac{L\,(1-r^2)-m\,(1+2r-r^2)}{(L-m)(3-r)(1-r^2)}
	\]}
	achieves the performance $\rate(\alg,F_{m,L}) = r$.
\end{theorem}

\begin{proof}
See \cref{app:proof_RAM_FmL}.
\end{proof}

\begin{remark}
	When $r$ is set to its minimum value of $1-\sqrt{\frac{m}{L}}$, the Robust Accelerated Method in \cref{thm:RAM} reduces to the Triple Momentum Method (see \cref{table:algos}).
\end{remark}

\subsection{Comparison with IQCs from robust control}
\label{sec:comparison}

We now compare our analysis in \cref{thm:lmi_cvx} for computing the worst-case convergence rate over $F_{m,L}$ with the use of integral quadratic constraints (IQCs)~\cite{iqc}.
Specifically, we apply the upper bound on convergence rate from the LMI in~\cite[Eq. 3.8]{lessard16} with the weighted off-by-one IQC in~\cite[Lemma 10]{lessard16}. The result is summarized below.

\begin{proposition}[Weighted off-by-one IQC]\label{prop:IQC}
Consider the algorithm $\alg = (A,B,C)$ defined in \eqref{eq:ssalg} satisfying \cref{ass:algform} applied to a function $f \in F_{m,L}$ defined in \cref{sec:intro}. If there exists $Q\succ 0$ satisfying
\[
    \bmat{\hat A^\tp Q \hat A - r^2Q & \hat A^\tp Q \hat B \\ \hat B^\tp Q \hat A & \hat B^\tp Q \hat B} + \bmat{\hat C & \hat D}^\tp \bmat{0 & 1 \\ 1 & 0} \bmat{\hat C & \hat D} \preceq 0,
\]
then $\rate(\alg,F_{m,L}) \leq r$, where the matrices $(\hat A,\hat B,\hat C,\hat D)$ are given by
\[
    \hat A = \bmat{A & 0 \\ -LC & 0}, \quad
    \hat B = \bmat{B \\ 1}, \quad
    \hat C = \bmat{LC & r^2 \\ -mC & 0}, \quad
    \hat D = \bmat{-1 \\ 1}.
\]
\end{proposition}
If $Q\succ 0$ satisfies the weighted off-by-one IQC LMI from \cref{prop:IQC}, then we can construct a feasible point for our analysis LMI from~\eqref{eq:lmi_cvx_rate} with $\ell=1$ as follows.
\begin{subequations}\label{eq:IQC_parameterization}
\begin{gather}
  p=2\,(L-m), \quad
  \Lambda_1 = \bmat{0 & 0 & 0 \\ r^2 & 0 & 0 \\ 1-r^2 & 0 & 0}, \quad
  \Lambda_2 = \bmat{0 & 0 & 0 \\ 0 & 0 & 1 \\ 0 & 0 & 0},\quad\text{and}\\
  P = \bmat{A & B \\ -LC & 1}^\tp Q \bmat{A & B \\ -LC & 1} - \bmat{LmC^\tp C & -mC^\tp \\ -mC & 1}.
\end{gather}
\end{subequations}
In this case, a Lyapunov function that certifies the linear rate $r$ is
\begin{equation}\label{eq:IQC_lyap}
  V(\xi^t,u^t,f^t) = \bmat{\xi^{t+1} \\ \zeta^{t+1}}^\tp\! Q \bmat{\xi^{t+1} \\ \zeta^{t+1}} + 2\,(L-m)\left(f^t - \frac{m}{2} \normm{y^t}^2\right) - \normm{u^t - m y^t}^2,
\end{equation}
where $\xi^{t+1} = A\xi^t + Bu^t$ and $\zeta^{t+1} \defeq u^t-L y^t$. Therefore, using the weighted off-by-one IQC can be interpreted as searching over this restricted class of Lyapunov functions.
Even though our analysis for computing the convergence rate in \cref{thm:lmi_cvx} is more general, the weighted off-by-one IQC formulation appears to be general enough to prove tight results. For example, while RAM was designed using the more general analysis, its Lyapunov function has the special form~\eqref{eq:IQC_lyap}.

In~\cite{scherer}, Zames--Falb multipliers~\cite{zames-falb} (of which the weighted off-by-one IQC is a special case) are used to formulate LMIs for computing both the convergence rate and the sensitivity. Just as with the weighted off-by-one IQC, using general Zames--Falb multipliers can also be interpreted as searching over a restricted class of Lyapunov functions, although a detailed comparison is beyond the scope of this work.

While the weighted off-by-one IQC formulation appears to achieve tight bounds on the convergence rate, computing tight bounds on the sensitivity requires the more general LMI~\eqref{eq:lmi_cvx_sensitivity}; see \cref{sec:computational_considerations} for further discussion.

\section{Numerical validation}\label{sec:numerical}

We perform numerical experiments to verify that our designs
(i) effectively trade off convergence rate and noise sensitivity,
(ii) use an adequate number of parameters (they are neither under- nor over-parameterized), and
(iii) outperform popular iterative schemes when applied to a worst-case test function.
Our code is available at \url{https://github.com/QCGroup/optalg}.

\subsection{Empirical validation}
\label{sec:pareto_validation}

To empirically validate our designs from \cref{thm:RHB,thm:RAM}, we performed a brute-force analysis of algorithms $(\alpha,\beta,\eta)$ and made a scatter plot of sensitivity vs.\ convergence rate. The following result facilitates sampling by providing bounds on admissible $(\alpha,\beta,\eta)$.

\begin{lemma}[parameter restriction]\label{cor:param_restriction}
	Consider a three-parameter algorithm $\alg=(\alpha,\beta,\eta)$ defined in \cref{eq:algform0}. Let $\F \in \{Q_{m,L},F_{m,L}\}$. If ${\rate(\alg,\F) < 1}$, then:
\begin{align*}
    0 &< \alpha < \frac{4}{L},
    &
    \frac{-2}{L-m} &< \alpha\eta < \frac{2}{L-m},
    &
    &\begin{cases}
        -1+L(\alpha\eta) < \beta < 1+m(\alpha\eta) &\text{if }\alpha\eta \geq 0 \\
        -1+m(\alpha\eta) < \beta < 1+L(\alpha\eta) & \text{if }\alpha\eta < 0.
    \end{cases}
\end{align*}
\end{lemma}

\begin{proof}
    See \cref{app:proof_param_bounds}.
\end{proof}

From \cref{cor:param_restriction}, we see that $\alpha$, $\alpha\eta$, and $\beta$ each have finite ranges. A convenient way to grid the space of possible $(\alpha,\beta,\eta)$ values is to first grid over $\alpha$, then $\alpha\eta$, then $\beta$, in a nested fashion, extracting associated $(\alpha,\beta,\eta)$ values at each step. Due to the multiplicative nature of the parameter $\alpha$, we opted to sample $\alpha$ logarithmically in the range $\bigl[10^{-5}, \tfrac{4}{L}\bigr]$, but to sample $\alpha\eta$ and $\beta$ linearly in their associated intervals.

\parag{Strongly convex quadratics (\texorpdfstring{$Q_{m,L}$}{Qml})}

\begin{figure}[ht]
	\centering
	\includegraphics[width=\linewidth]{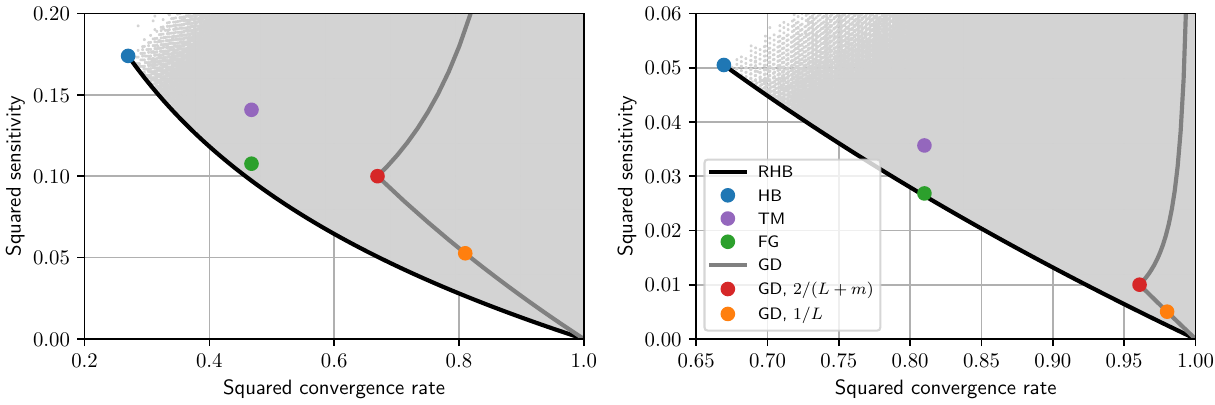}
    
\caption{Plot of sensitivity squared vs.\ rate squared for algorithms applied to the function class $Q_{1,10}$ (left panel) and $Q_{1,100}$ (right panel), found using \cref{cor:perf_quad}. Each point in the point cloud corresponds to an algorithm $(\alpha,\beta,\eta)$.
The Pareto-optimal front coincides with our proposed Robust Heavy Ball method (RHB, \cref{thm:RHB}), tuned using $r \in \bigl[ \frac{\sqrt{L}-\sqrt{m}}{\sqrt{L}+\sqrt{m}}, 1\bigr]$ to mediate the trade-off.}
\label{fig:QmL_cloud}
\end{figure}

We show our brute-force search for the class $Q_{m,L}$ in \cref{fig:QmL_cloud} for $Q_{1,10}$ and $Q_{1,100}$. For this figure, we used the sampling approach based on \cref{cor:param_restriction} with $500\times 201\times 200$ samples for $\alpha$, $\alpha\eta$, and $\beta$, respectively.

In \cref{fig:QmL_cloud}, each algorithm $(\alpha,\beta,\eta)$ corresponds to a single gray dot
The RHB curve shows each possible tuning as we vary the parameter $r$. We observe that RHB perfectly traces out the boundary of the point cloud, which represents the Pareto-optimal algorithms. In other words, for a fixed convergence rate $r$, RHB with parameter $r$ achieves this rate and is also as robust as possible to additive gradient noise.
\looseness=-1

\cref{fig:QmL_cloud} reveals that GD with $0 < \alpha < \tfrac{2}{L+m}$ is outperformed by RHB on the function class $Q_{m,L}$. We also plot the performance of GD for $\alpha > \tfrac{2}{L+m}$, which is even worse as this leads to slower convergence \emph{and} increased sensitivity. \cref{fig:QmL_cloud} also reveals that FG is strictly suboptimal compared to RHB based on the analysis in \cref{thm:lmi_cvx}, although the optimality gap appears to shrink as $L/m$ gets larger.

\parag{Smooth strongly convex functions (\texorpdfstring{$F_{m,L}$}{Fml})}

\begin{figure}[htb]
	\centering
	\includegraphics[width=\linewidth]{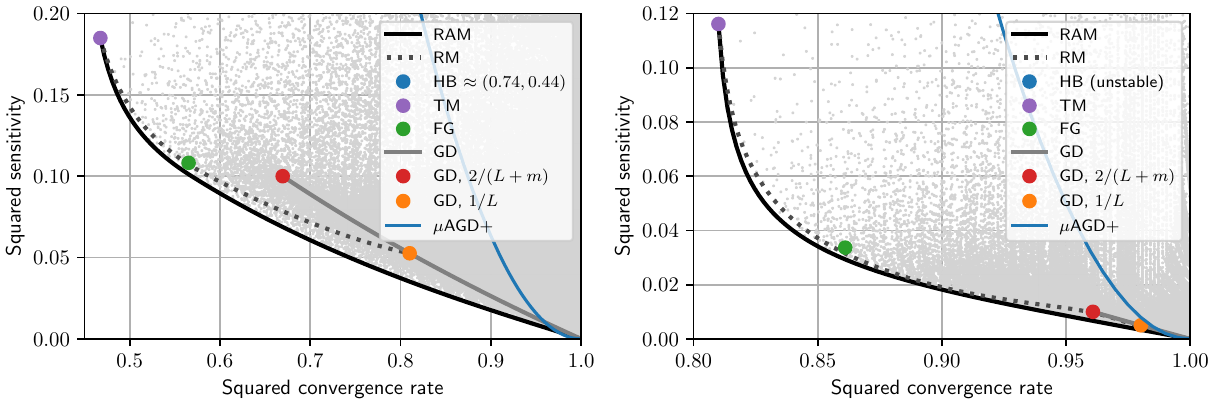}
    
\caption{Plot of squared sensitivity vs.\ squared rate for algorithms applied to the function class $F_{1,10}$ (left panel) and $F_{1,100}$ (right panel), found using \cref{thm:lmi_cvx}. Each point in the point cloud corresponds to an algorithm $(\alpha,\beta,\eta)$.
RM and GD are strictly suboptimal, and the Pareto-optimal front closely matches RAM (\cref{thm:RAM}), tuned using $r \in \bigl[ 1-\sqrt{\frac{m}{L}}, 1\bigr]$ to mediate the trade-off. For reference, the bound \eqref{eq:AGD} on $\mu$AGD+ from~\cite{cohen2018} is also shown.}
\label{fig:FmL_cloud}
\end{figure}

We show our brute-force search for the class $F_{m,L}$ in \cref{fig:FmL_cloud} for $F_{1,10}$ and $F_{1,100}$. For this figure, we used the same sampling approach as in \cref{fig:QmL_cloud}, with $200\times51\times50$ samples. When applying \cref{thm:lmi_cvx}, we used a lifting dimension $\ell=1$ to compute the rate and $\ell=6$ to compute the sensitivity. For more details on these choices, see \cref{sec:computational_considerations}.

The Robust Momentum (RM) method from~\cite{rmm} interpolates between TM and GD with $\alpha=\frac{1}{L}$ and does trade off convergence rate and sensitivity, but based on our analysis in \cref{thm:lmi_cvx} it is strictly outperformed by our proposed RAM. The gap in performance between RM and RAM appears to shrink as $L/m$ gets larger.

While \cref{fig:FmL_cloud} indicates that RAM effectively trades off rate and robustness, it is strictly suboptimal\footnote{Suboptimality is difficult to verify for the function class $F_{m,L}$ since the results depend on the lifting dimension $\ell$. However, we performed extensive numerical computations to ensure that $\ell$ is sufficiently large; see \cref{sec:computational_considerations}.}. Suboptimality becomes most apparent when $L/m$ is small and $r$ is close to $1$. For example, consider RAM with the parameter choice $r=0.9$, $m=1$, and $L=2$, which corresponds to $(\alpha,\beta,\eta) = (0.019, 0.66, -3.631579)$. Solving the LMIs in \cref{thm:lmi_cvx} yields the rate $0.9000$ and sensitivity $0.22057$. However, if we change $\eta$ and use the tuning $(\alpha,\beta,\eta) = (0.019, 0.66, 0.00)$ instead, we obtain the same rate with the strictly smaller sensitivity $0.1676$. Larger optimality gaps can be found by making $L/m$ even closer to $1$, although such cases are not practical.

\begin{remark}
The left panel of \cref{fig:FmL_cloud} ($L/m=10$) shows a denser point cloud than the right panel ($L/m=100$), despite using the same number of sample points. This occurs because the right panel spans a larger range of $\gamma$ values (the vertical axis is truncated), indicating that desirable algorithm tunings become harder to find by random sampling as $L/m$ increases.
\end{remark}

\subsection{Justifying the three-parameter algorithm family}\label{sec:justification_three_param}

\begin{figure}[htb]
	\centering
	\includegraphics[width=\linewidth]{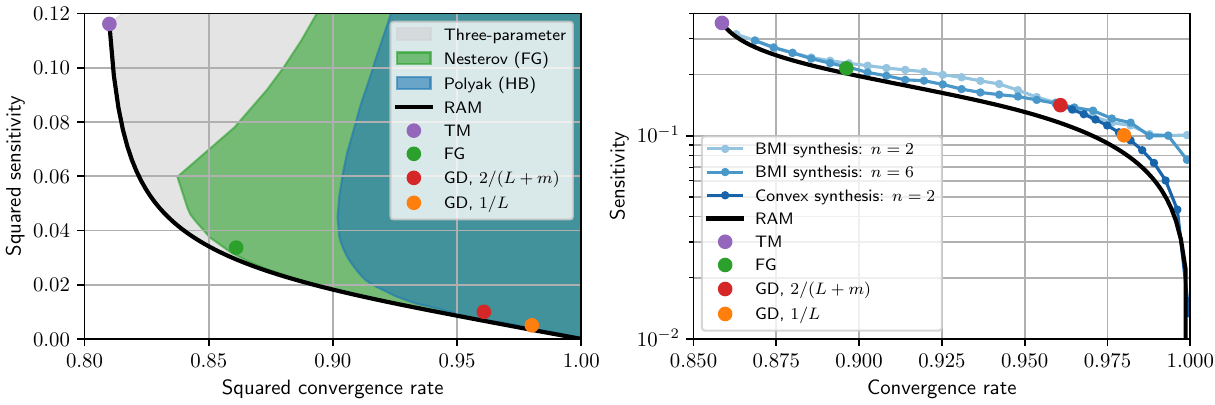}
    
\caption{\textbf{Left:} Regions of the sensitivity vs.\ rate trade-off space for $F_{1,100}$ covered by the three-parameter family $(\alpha,\beta,\eta)$, the Nesterov (Fast Gradient) family $(\alpha,\beta,\beta)$, and the Polyak (Heavy Ball) family $(\alpha,\beta,0)$. The FG and HB families are not expressive enough to capture the whole trade-off space. \textbf{Right:} Comparison of RAM with the numerically synthesized algorithms (using a state dimension up to $n=6$) from~\cite{scherer} for $F_{1,50}$. RAM outperforms despite using only two states of memory. We plot the log of the normalized sensitivity vs. the rate to match \cite[Fig.~6]{scherer}.}
\label{fig:venn_scherer}
\end{figure}

A natural question to ask is whether something as general as our three-parameter family~\eqref{eq:algform0} is needed to achieve the performance of our designs. Several recent works have restricted their attention to optimizing algorithms with two parameters $(\alpha,\beta)$ in either Nesterov's FG or Polyak's HB form \cite{mert,mihailo,HB:convex}.
From our results in \cref{sec:algodesign_QmL,sec:pareto_validation}, the HB form is sufficient for the class $Q_{m,L}$. However, neither the HB or the FG forms are sufficient for $F_{m,L}$. While some algorithms in these restricted classes achieve acceleration, they are incapable of obtaining the same trade-off between convergence rate and sensitivity as a properly-tuned three-parameter method, as illustrated in \cref{fig:venn_scherer} (left panel). Indeed, even when there is no noise, no algorithm in the FG or HB families achieves the optimal convergence rate for the function class $F_{m,L}$, which is attained by Van Scoy et al.'s Triple Momentum (TM) method~\cite{tmm}.

Alternatively, we could ask whether using more than three parameters could lead to further improvement. As explained in~\cref{sec:parameterization}, any algorithm with $n=2$ states can be represented by three parameters. In general, we would need $2n-1$ parameters to represent an algorithm with $n$ states. In principle, our methodology of \cref{sec:design} can still be applied, but the associated semidefinite programs become substantially more difficult to solve and we were unable to find better designs.

An alternative approach, presented in \cite{scherer}, used a convex synthesis procedure and bilinear matrix inequalities to numerically construct algorithms that trade off convergence rate and sensitivity. As shown in \cref{fig:venn_scherer} (right panel), these synthesized algorithms do not outperform RAM, despite using up to $n=6$ states.

\subsection{Simulation of a worst-case test function}
\label{sec:simulation_worst-case}

\begin{figure}[htb]
	\centering
	\includegraphics[width=\linewidth]{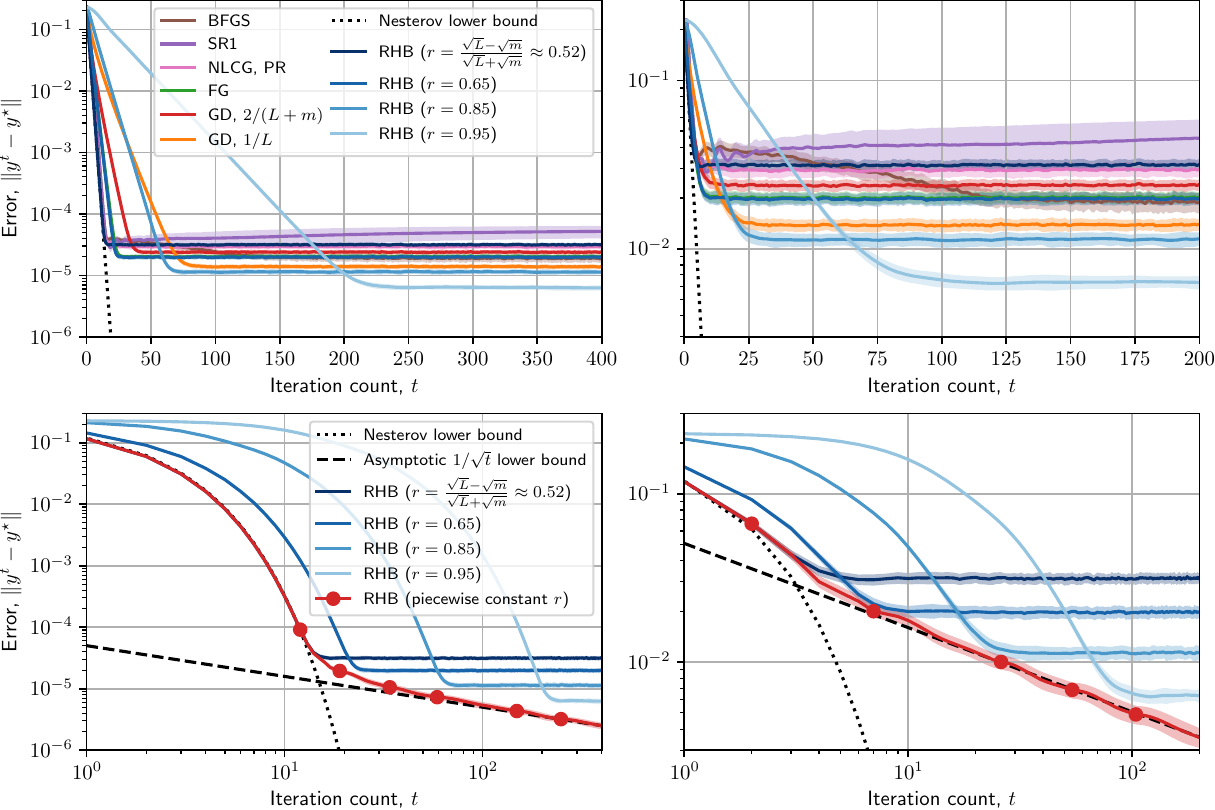}
    
  \caption{Simulations of various algorithms with low noise ($\sigma=10^{-5}$, left column) and higher noise ($\sigma=10^{-2}$, right column). Each algorithm was simulated on Nesterov's lower-bound quadratic function, with $m=1$, $L=10$, and dimension $d=100$. Shaded regions indicate $\pm 1$ standard deviations about the mean across $100$ trials (different noise realizations). Different tunings of our proposed Robust Heavy Ball (RHB) from \cref{thm:RHB} effectively trade off convergence rate and steady-state error (sensitivity to noise). Bottom row: the red curve shows RHB with (hand-tuned) piecewise constant $r$, where the red dots indicate switch points.}
  \label{fig:simulation_quadratics}
\end{figure}

We simulated various algorithms on Nesterov's lower-bound function, which is a quadratic with a tridiagonal Hessian~\cite[\S2.1.4]{nesterov_book}. We used $d=100$ with $m=1$ and $L=10$ and initialized each algorithm at zero. The results are reported in \cref{fig:simulation_quadratics}. We tested both a \emph{low noise} ($\sigma=10^{-5}$, left column) and a \emph{higher noise} ($\sigma=10^{-2}$, right column) regime. We recorded the mean and standard deviation of the error across $100$ trials for each algorithm (the trials differ only in the noise realization). 

\cref{fig:simulation_quadratics} shows that RHB from \cref{thm:RHB} trades off convergence rate (the initial decreasing slope) with sensitivity to noise (the steady-state value), and compares favorably to other methods. The other methods we tested (first row of \cref{fig:simulation_quadratics}) are generally suboptimal compared to RHB, in the sense that there is some choice of tuning parameter $r$ such that RHB is both faster and has smaller steady-state error.

In addition to GD and FG, we tested Nonlinear Conjugate Gradient (NLCG) with the Polak-Ribi\`ere (PR) update scheme.\footnote{We also tested other popular NLCG update schemes such as: Fletcher--Reeves, Hestenes--Stiefel, and Dai--Yuan; all produced similar trajectories to PR.} NLCG performs similarly to RHB with the most aggressive tuning, which is equivalent to the Heavy Ball method. 
We also tested the popular quasi-Newton methods~\cite{nocedal_wright}
Broyden--Fletcher--Goldfarb--Shanno (BFGS) and Symmetric Rank-One (SR1), which performed strictly worse than RHB.\footnote{
Both NLCG and BFGS use line search; given a current point $y \in \R^d$ and search direction $s \in \R^d$, they search for $\alpha \in \R$ that minimizes $f(y + \alpha s)$.
In practice, \emph{inexact} line searches are performed at each timestep with a stopping criterion such as the Wolfe conditions. To show these algorithms in the most charitable light, we used exact line searches but substituted the noisy gradient oracle.
Specifically, with $f(y) = \frac{1}{2}(y-y^\star)^\tp Q (y-y^\star)$, the optimal stepsize is $\alpha^\star = -(s^\tp \df(y))/(s^\tp Q s)$. We used this formula, but replaced $\df(y)$ by the noisy  $\df(y) + w$ (exact knowledge of $Q$ but not $y^\star$).}

In the second row of \cref{fig:simulation_quadratics}, we use the same settings as the first row, except iterations are plotted on a log scale as in \cref{fig:demo_rand_trajectories}. We also show a hand-tuned version of RHB with piecewise constant $r$.\footnote{We use a hand-tuned schedule to illustrate our results; more systematic scheduling methods have been proposed, such as the \textsc{restart+slowdown} method in \cite{cohen2018}.} Whenever $r$ is changed, we re-initialize the algorithm via $x^{t-1}=x^t$.
In the transient phase, the hand-tuned RHB matches Nesterov's lower bound \eqref{eq:nesterov_lowerbound}.
In the stationary phase, it matches the asymptotic lower bound (slope of $-1/2$) described in \cref{sec:background}.\footnote{The lower bound is $1/t$ for the squared error, hence $1/\sqrt{t}$ for the error, which appears as a line of slope $-1/2$ on the log-log scale of \cref{fig:simulation_quadratics}, bottom row.}

\subsection{Computational considerations}\label{sec:computational_considerations}

We used Julia v.1.12.1~\cite{julia} for all computation, and our code is available at
\url{https://github.com/QCGroup/optalg}.
Applying \cref{thm:lmi_cvx} required picking $\ell$. Using larger $\ell$ could reduce conservatism, but at the cost of larger LMIs.
To this end, we performed a pilot study using the arbitrary-precision solver Hypatia~\cite{hypatia}.
In \cref{table:dimverif}, we apply \cref{thm:lmi_cvx} and report the performance of FG with standard tuning (see \cref{table:algos}), with 30 significant digits\footnote{For each case, we computed optimal primal and dual solutions, verified that each was strictly feasible, and ensured the duality gap was less than $10^{-30}$.}.
The rate did not improve beyond $\ell=1$ but the sensitivity bound continued to improve.

\begin{table}[ht]
	\centering
	\caption{Performance bounds for FG with standard tuning (see \cref{table:algos}) on $F_{1,100}$ found by solving the LMIs from \cref{thm:lmi_cvx} using different lifting dimensions $\ell$ and arbitrary-precision arithmetic.
    All digits shown are significant, and digits that represent conservatism due to choice of $\ell$ are in red.
    The last digit is rounded up so all numbers shown represent valid upper bounds.}
  \label{table:dimverif}
    \newcommand{\bad}[1]{\textcolor{red}{#1}}
	\renewcommand\arraystretch{1}
    \setlength{\tabcolsep}{5pt}
    { \begin{tabular}[t]{lll}
		\toprule
		$\ell$ & Rate upper bound & Sensitivity upper bound\\ 
		\midrule
%
%
%
        1 & 0.927933110963822850550405971295 & 0.\bad{200765554941665409665170222653} \\
        2 & 0.927933110963822850550405971295 & 0.18\bad{5908689372372655449415199727} \\
        3 & 0.927933110963822850550405971295 & 0.183\bad{728741557681344670475184476} \\
        4 & 0.927933110963822850550405971295 & 0.183\bad{513628727818276451861705856} \\
        5 & 0.927933110963822850550405971295 & 0.18349\bad{3246027954433809754833544} \\
        6 & 0.927933110963822850550405971295 & 0.183491\bad{309055074730945919778091} \\
        7 & 0.927933110963822850550405971295 & 0.1834911\bad{24755530197869873815523} \\
        8 & 0.927933110963822850550405971295 & 0.18349110\bad{7215459968837777689264} \\
        9 & 0.927933110963822850550405971295 & 0.183491105\bad{546079662576734624809} \\
        10& 0.927933110963822850550405971295 & 0.1834911053\bad{87195122067178179096} \\
		\bottomrule
  \end{tabular}}
\end{table}

To generate \cref{fig:FmL_cloud}, we solved the LMIs from \cref{thm:lmi_cvx} using JuMP v.1.29.2 \cite{JUMP} with the Mosek solver v.11.0.1 \cite{mosek} and default settings. We used $\ell=1$ for the rate, with a bisection search tolerance of $10^{-6}$, and $\ell=6$ for the sensitivity. With these choices, finding the rate and sensitivity of a given three-parameter algorithm with 5 significant digits of precision each took about 100~ms on a conventional laptop.

\section{Concluding remarks}\label{sec:conclusion}

For $Q_{m,L}$ and $F_{m,L}$, we provided (i) efficient methods for computing the convergence rate and noise sensitivity for a broad class of first-order methods, and (ii) first-order algorithms designs, each with a single tunable parameter that directly trades off convergence rate and sensitivity to additive gradient noise.

An interesting future direction is exploring adaptive versions of these algorithms, where the parameter $r$ is varied over time. We showed in \cref{fig:simulation_quadratics} that a hand-tuned piecewise constant version of RHB can match both Nesterov's lower bound and the gradient lower bound, so more sophisticated adaptive schemes such as those described in \cref{sec:background} might also work.
One could also adjust parameters continually, but proving the convergence of adaptive algorithms is generally more challenging.
For example, the well-known ADMM algorithm is often tuned adaptively to improve transient performance, even when convergence guarantees only hold for fixed parameters~\cite[\S3.4.1.]{boyd_admm}.
Nevertheless, LMI-based approaches have been successfully used to prove convergence of algorithms with time-varying parameters~\cite{diss_ICML,fazlyab}.

Another interesting open question is whether our analysis is tight. For  $F_{m,L}$ (\cref{thm:lmi_cvx}), our bounds depend on $\ell$. It is unknown (i) how large $\ell$ needs to be in order to obtain the tightest possible bounds on convergence rate and sensitivity, and (ii) whether \cref{thm:lmi_cvx} always produces tight bounds as $\ell\to\infty$.

\appendix

\section{Detailed proofs}

\subsection{Proof of \texorpdfstring{\cref{lem:alg3} (3-parameter family)}{3-parameter family}}\label{app:proof_3_param}
The three-parameter family~\eqref{eq:algform0} is of the general form~\eqref{eq:ssalg} with $n=2$;
\begin{align}\label{eq:algform}
\hat\xi^t &= \bmat{x^t \\ x^{t-1}}, &
\hat A &= \bmat{1+\beta & -\beta \\ 1 & 0}, &
\hat B &= \bmat{-\alpha \\ 0}, &
\hat C &= \bmat{1+\eta & -\eta}.
\end{align}
Let $(A,B,C)$ be a general algorithm of the form \eqref{eq:ssalg} satisfying \cref{ass:algform} with $n=2$.
We will construct a transformation matrix $T$, as introduced in \cref{rem:state_transformation}, that transforms $(A,B,C)$ to the form of \eqref{eq:algform}.

Let the eigenvalues of $A$ be $\{1,\beta\}$ with corresponding eigenvectors $v_1,v_2\in\R^2$. First transform using $T_1 = \bmat{v_1 & v_2}^{-1}$, which has the effect of diagonalizing $A$, and obtain
$(A,B,C) \mapsto \left( \sbmat{1 & 0 \\ 0 & \beta}, \sbmat{b_1 \\ b_2}, \sbmat{c_1 & c_2}\right)$. Then, transform using $T_2 = \frac{c_1}{b_2}\sbmat{b_2 & -\beta b_1\\b_2 & -b_1}$ and obtain\\
$\left( \sbmat{1+\beta & -\beta \\ 1 & 0}, \sbmat{ (1-\beta)b_1c_1\\0}, \sbmat{1+\frac{\beta b_1c_1+b_2c_2}{(1-\beta)b_1c_1} & -\frac{\beta b_1c_1+b_2c_2}{(1-\beta)b_1c_1}}\right)$. Setting $\alpha = -(1-\beta)b_1c_1$ and $\eta = \frac{\beta b_1c_1+b_2c_2}{(1-\beta)b_1c_1}$, the transformed system matches \eqref{eq:algform}. Thus, applying the transformation $T = T_2T_1$ transforms a generic $(A,B,C)$ to the form \eqref{eq:algform0} and if we initialize the new system with state $\hat \xi_0 = T \xi_0$, where $\xi_0$ is the initial state of $(A,B,C)$, both systems will produce the same sequence of iterates $(y^0,y^1,\dots)$.\qed

\subsection{Proof of \texorpdfstring{\cref{prop:perf_quad_general} ($Q_{m,L}$ analysis)}{quadratic analysis}}\label{app:proof_QmL}

Consider a particular $f\in Q_{m,L}$.
We can write $\df(y) = (y-y^\star) Q$ for some $Q \in \R^{d\times d}$ satisfying $Q = Q^\tp \succ 0$ (recall that the iterates are row vectors, so $y\in \R^{1\times d}$). The dynamics~\eqref{eq:ssalg} become
\begin{align}\label{eq:dynamics-quadratic}
	\xi^{t+1}-\xi^\star = A(\xi^t-\xi^\star) + BC(\xi^t-\xi^\star) Q + B w^t
	\quad\text{and}\quad
	y^t-y^\star = C (\xi^t-\xi^\star).
\end{align}
If we diagonalize $Q$ as $V \Lambda V^\tp$ with $V = \bmat{v_1 & \ldots & v_d}$ and $\Lambda = \diag(q_1,\ldots,q_d)$, we can split the dynamics into $d$ decoupled systems by multiplying on the right by $v_i$ and letting $\hat \xi^t_i \defeq (\xi^t-\xi^\star) v_i$ and similarly for $\hat y_i^t$ and $\hat w_i^t$, obtaining
\begin{align}\label{eq:ss_decoupled}
	\hat \xi^{t+1}_i = (A+q_i B C)\hat \xi^t_i + B \hat w^t_i
  \qquad\text{and}\qquad
	\hat y^t_i = C \hat \xi^t_i.
\end{align}
We now have $\hat \xi^t_i \in \R^{n\times 1}$, $\hat y^t_i \in \R$, and $\hat w^t_i\in \R$.
The rate and sensitivity of \eqref{eq:dynamics-quadratic} can be found by analyzing the simpler system~\eqref{eq:ss_decoupled}. Specifically, the rate for $\alg$ applied to $f$ is $\max_{i\in\{1,\dots,d\}}\rho(A+q_i B C)$, so the worst case over $Q_{m,L}$ is the supremum of $\rho(A+qBC)$ over all $q\in[m,L]$ (the set of possible eigenvalues of $Q$), as required.

For sensitivity, define $X_i^t \defeq \E\bigl[ (\hat\xi_i^t)(\hat\xi_i^t)^\tp \bigr] \in\R^{n\times n}$,
$W_i^t \defeq \E\left[ (\hat w_i^t)^2 \right] \in \R$, and
$Y_i^t \defeq \E\left[ (\hat y_i^t)^2 \right] \in \R$.
By assumption, $\sum_{i=1}^d W_i^t \leq \sigma^2$ for all $t$ and $W_i^t \geq 0$ for all $i,t$.
The worst-case sensitivity for $\alg$ applied to $f$ is achieved when all weight is placed on a single index $i$ for all $t$, so $W_i^t = W_i = \sigma^2$. For this $i$, we can write the recurrence
\begin{equation}\label{qq1}
X_i^{t+1} = (A+q_iBC)X_i^t (A+q_iBC)^\tp + \sigma^2 B B^\tp
\quad\text{and}\quad
Y_i^t = C X_i^t C^\tp.
\end{equation}
Define the Lyapunov operator $\mathcal{L}_i(Z) \defeq \sum_{k=0}^\infty (A+q_i BC)^k Z \bigl((A+q_iBC)^\tp\bigr)^k$. This series converges whenever $\rho(A+q_i BC)<1$ to the (unique) solution $X$ of the Lyapunov equation $(A+q_iBC) X(A+q_iBC)^\tp-X +Z=0$. Since $\rate(\alg,Q_{m,L})<1$ by assumption, $\mathcal{L}_i$ is well-defined for all $q_i \in [m,L]$. Applying $\mathcal{L}_i$ above, we obtain
\[
\limsup_{T\to\infty} \E_{w\sim \mathbb{P}} \frac{1}{T\sigma^2}\sum_{t=0}^{T-1} \norm{y^t-y^\star}^2 
\leq \max_{1\leq i\leq d} \frac{1}{\sigma^2} C \mathcal{L}_i(\sigma^2 BB^\tp) C^\tp
= \max_{1\leq i\leq d} C \mathcal{L}_i(BB^\tp) C^\tp,
\]
where the final equality is due to linearity of $\mathcal{L}_i$.
The worst case over $Q_{m,L}$ is then the supremum of $\sqrt{C\mathcal{L}(BB^\tp)C^\tp}$ over all possible eigenvalues $q\in[m,L]$ of~$Q$. Define the adjoint Lyapunov operator $\mathcal{L}^*(Z) \defeq \sum_{k=0}^\infty \bigl((A+qBC)^\tp\bigr)^k Z (A+q BC)^k$ as the unique solution $Y$ to the Lyapunov equation $(A+qBC)^\tp Y(A+qBC)-Y +Z=0$. Applying Lagrangian duality, we have $\sqrt{C\mathcal{L}(BB^\tp)C^\tp} = \sqrt{B^\tp \mathcal{L}^*(C^\tp C) B}$.\qed

\subsection{Proof of \texorpdfstring{\cref{cor:perf_quad} (reduced analysis for $Q_{m,L}$)}{reduced quadratic analysis}}\label{app:proof_QmL_reduced}

  We first prove the rate expression~\eqref{eq:rho_mL_quad}. Substituting the algorithm form~\eqref{eq:algform} into \cref{prop:perf_quad_general},
  \begin{equation}\label{dolphin}
  \rate(\alg,Q_{m,L}) = \sup_{q \in[m,L]} \phi(q) \defeq \rho\left(
    \bmat{\beta +1-\alpha  (\eta +1) q \quad & -\beta+\alpha  \eta  q   \\
    1 & 0}\right).
  \end{equation}
  We will prove that $\phi(q)$ defined above is quasiconvex~\cite[\S3.4]{boyd_vandenberghe}. The characteristic polynomial of the matrix in~\eqref{dolphin} is $\chi(z) = 
  z^2+(\alpha  (\eta +1) q  -\beta -1)z +(\beta-\alpha  \eta  q )$. Given a polynomial with real coefficients, a necessary and sufficient condition for its roots to lie inside the unit circle is given by the Jury test~\cite[\S4.5]{fadali2013digital}. For the quadratic $z^2 + a_1 z + a_0$, the Jury test reduces to $1 + a_1 + a_0 > 0$, $1 - a_1 + a_0 > 0$, and $-1 < a_0 < 1$. Applying the Jury test to $\chi(rz)$, we find that $\phi(q) < r$ if and only if\looseness=-1
  \begin{equation}\label{eq:jury_rho}
  \begin{aligned}
    (1-r) (\beta - r ) + \alpha (\eta  r - \eta + r ) q  &> 0, &
    r^2+\beta -\alpha \eta q &> 0, \\
    (1+r) (\beta + r ) - \alpha (\eta  r + \eta + r ) q  &> 0, &
    r^2-\beta + \alpha \eta q &> 0.
  \end{aligned}
  \end{equation}
  The inequalities~\eqref{eq:jury_rho} are linear in $q$, so the sublevel sets $\{q\;\vert\;\phi(q) < r\}$ are open intervals. Therefore, $\phi$ is quasiconvex and attains its supremum over $q \in [m,L]$ at one of the endpoints, $q=m$ or $q=L$. The explicit formula for $\phi(q)$ can be found by applying the quadratic formula to find the roots of $\chi(z)$. In~\eqref{eq:rho_mL_quad}, $\Delta$ is the discriminant of $\chi(z)$ and the two cases correspond to whether the roots are real or complex.

  We next prove the expression for the sensitivity in~\eqref{eq:gamma_mL_quad}. Substituting the algorithm form~\eqref{eq:algform} into \cref{prop:perf_quad_general}, we can explicitly solve the linear equation for $P_q $ in \eqref{eq:gamma_mL_quad} to obtain $h(q) = B^\tp P_q B$, where $h(q)$ is defined in~\eqref{eq:gamma_mL_quad}. When $r=1$, the Jury conditions~\eqref{eq:jury_rho} reduce to:
  \begin{subequations}\label{eq:jury_one}
    \begin{align}
      \alpha q  &> 0, \label{ea}\\
      2\beta +2 -\alpha (2\eta +1 ) q  &> 0, \label{eb}\\
      1+\beta -\alpha  \eta  q   &>0, \label{ec}\\
      1-\beta + \alpha  \eta  q  &>0. \label{ed}
    \end{align}
  \end{subequations}
When the rate is strictly less than one, $h(q)$ is positive and convex. To see why, evaluate $h(q)$ and $h''(q)$. After routine algebraic manipulations, we obtain:
  \begin{subequations}\label{turtle}
  \begin{align}
    h(q) &= \frac{\alpha ^2 (2 \eta +1)^2}{2 (1-\beta + \alpha  \eta  q ) (2\beta +2 -\alpha (2\eta +1 ) q )}
    +\frac{\alpha^2 }{2 \alpha q   (\alpha  \eta  q  -\beta +1)},\label{turtle1}\\
    h''(q) &=
    \frac{\alpha ^4 (2 \eta +1)^2 \Bigl(3 (2 \alpha  \eta  (2 \eta +1) q  -4 \beta  \eta -\beta +1)^2+(4 \eta-\beta+1)^2\Bigr)}{4 (1-\beta + \alpha  \eta  q )^3 (2\beta +2 -\alpha (2\eta +1 ) q )^3} \notag\\ 
    &\hspace{5cm}+
    \frac{\alpha q \Bigl(3 \bigl(2\alpha  \eta  q  +1-\beta\bigr)^2+\bigl(1-\beta \bigr)^2\Bigr)}{4 q  ^4 (\alpha  \eta q -\beta +1)^3}. \label{turtle2}
  \end{align}
  \end{subequations}
  In the form~\eqref{turtle}, it is clear that whenever the rate is strictly less than one (i.e., when \eqref{eq:jury_one} holds) we have $h(q) > 0$ and $h''(q)> 0$. Therefore, the quantity under the square root in~\eqref{eq:gamma_mL_quad} is always positive and $h(q)$ is convex, so it attains its supremum over $q\in[m,L]$ at one of the endpoints, $q=m$ or $q=L$.\qed

\subsection{Proof of \texorpdfstring{\cref{thm:lmi_cvx} ($F_{m,L}$ analysis)}{analysis result for smooth strongly convex case}}\label{app:proof_FmL_analysis}

  Consider a trajectory of the dynamics $(\xi^t,u^t,y^t,w^t)$ \eqref{eq:ssalg} with $w^t=0$.
  Form the augmented vectors and state as described in \cref{sec:perf_cvx}. Multiply the LMIs \eqref{eq:lmi_cvx_1} and \eqref{eq:lmi_cvx_3} on the right and left by $\col(\bxi^t,u^t)\in\R^{(n+2\ell)\times d}$ and its transpose, respectively, and take the trace. Also, take the inner product of \eqref{eq:lmi_cvx_2} and \eqref{eq:lmi_cvx_4} with $\f^t$, which preserves the inequality because $\f^t$ is elementwise nonnegative. The resulting inequalities are:
\begin{subequations}\label{eq:rate_cvx}
	\begin{align}\label{eq:rate_cvx_a}
		\trace\,(\bxi^{t+1})^\tp P \bxi^{t+1} - r^2 \trace\,(\bxi^t)^\tp P \bxi^t + \trace \bmat{\y^t \\ \u^t}^\tp \Pi(\Lambda_1) \bmat{\y^t \\ \u^t} &\leq 0, \\
    \label{eq:rate_cvx_b}
		p^\tp (Z_+ - r^2 Z) \f^t + \pi(\Lambda_1)^\tp \f^t &\leq 0, \\
    \label{eq:rate_cvx_c}
		\|\xi^t-\xi^\star\|^2 - \trace\,(\bxi^t)^\tp P \bxi^t + \trace \bmat{\y^t \\ \u^t}^\tp \Pi(\Lambda_1) \bmat{\y^t \\ \u^t} &\leq 0, \\
		\label{eq:rate_cvx_d}
		-p^\tp Z \f^t + \pi(\Lambda_2)^\tp \f^t &\leq 0.
	\end{align}
\end{subequations}
Summing \eqref{eq:rate_cvx_a}$+$\eqref{eq:rate_cvx_b} and \eqref{eq:rate_cvx_c}$+$\eqref{eq:rate_cvx_d} and applying~\eqref{eq:cvx_ineq}, we recover the lower bound and decrease properties in \cref{lem:lyap_rate} and the rate bound result follows.

For the second part of the proof, we do not restrict $w^t=0$, and perform similar operations to the inequalities \eqref{eq:lmi_cvx_sensitivity} as in the first part to obtain the inequalities
\begin{subequations}\label{eq:sensitivity_cvx}
	\begin{align}
		\trace\,(\bxi^{t+1}-\H w^t)^\tp P (\bxi^{t+1}-\H w^t) - \trace\,(\bxi^t)^\tp P \bxi^t + \trace \bmat{\y^t \\ \u^t}^\tp \Pi(\Lambda_1) \bmat{\y^t \\ \u^t} + \|\tilde y^t\|^2 &\leq 0, \label{eq:sensitivity_cvx_a} \\
		p^\tp (Z_+ - Z) \f^t + \pi(\Lambda_1)^\tp \f^t &\leq 0,
    \label{eq:sensitivity_cvx_b} \\
		-\trace\,(\bxi^t)^\tp P \bxi^t + \trace \bmat{\y^t \\ \u^t}^\tp \Pi(\Lambda_2) \bmat{\y^t \\ \u^t} &\leq 0,
    \label{eq:sensitivity_cvx_c} \\
		-p^\tp Z \f^t + \pi(\Lambda_2)^\tp \f^t &\leq 0,
    \label{eq:sensitivity_cvx_d}
	\end{align}
\end{subequations}
where $\tilde y^t = y^t-y^\star$. Summing \eqref{eq:sensitivity_cvx_a}$+$\eqref{eq:sensitivity_cvx_b} and \eqref{eq:sensitivity_cvx_c}$+$\eqref{eq:sensitivity_cvx_d}, applying~\eqref{eq:cvx_ineq}, and using the definition of the augmented state~\eqref{eq:state_aug} and the Lyapunov function~\eqref{eq:lyap},
\begin{align*}
	V(\x^{t+1}) - V(\x^t) + \|\tilde y^t\|^2 - 2\,\trace\,(\A\bxi^t+\B u^t)^\tp P \H w^t - \trace (w^t)^\tp \H^\tp P \H w^t &\leq 0
\end{align*}
and $V(\x^t) \geq 0$. Taking the expectation of both inequalities, the term $-2(\A\bxi^t+\B u^t)^\tp P \H w^t$ in the first inequality vanishes because $w^t$ is zero-mean and is independent of $\bxi^t$ and $u^t$, which only depend on $w^{t-1},w^{t-2},\dots$. Also, since $w^t$ has expected squared norm bounded by $\sigma^2$, we have $\trace\E\bigl( (w^t)^\tp \H^\tp P \H w^t\bigr) \leq \sigma^2 (\H^\tp P \H)$. Thus, the previous inequalities imply the decrease and lower bound conditions in \cref{lem:lyap_sensitivity} with $\gamma^2 = \sigma^2 (\H^\tp P\H)$, which implies the bound on the sensitivity.\qed

\subsection{Proof of \texorpdfstring{\cref{thm:RGM2} (GD analysis for $F_{m,L}$)}{GD analysis for smooth strongly convex case}}\label{app:proof_GD_FmL}

Setting $\ell=0$ in \cref{thm:lmi_cvx} and applying \cref{rem:ellzero} with the GD parameters from \eqref{eq:ssalg} ($A=C=1$ and $B=-\alpha$), the LMIs~\eqref{eq:lmi_cvx_rate}--\eqref{eq:lmi_cvx_sensitivity} with stepsize $\alpha=\frac{1}{m}(1-r)$ are satisfied with $P_1=P_2=\frac{1}{1-r^2}$ and $\lambda_1=\lambda_2=\frac{r}{m(L-m)(r +1)}$ for all $\frac{L-m}{L+m} \le r < 1$. Therefore, $\rate(\alg,F_{m,L}) \leq r$ and $\sensitivity(\alg,F_{m,L}) \leq \alpha \sqrt{P_2} = \frac{1}{m}\sqrt{\frac{1-r}{1+r}}$. Since these upper bounds match the exact values found for $Q_{m,L}$ in \cref{thm:RGM}, we conclude the upper bounds are tight, with the worst-case function being the quadratic $f(y) = \tfrac{m}{2}\norm{y}^2$.
\qed

\subsection{Proof of \texorpdfstring{\cref{thm:RAM} (RAM analysis for $F_{m,L})$}{RAM analysis for smooth strongly convex case}}\label{app:proof_RAM_FmL}

To prove that RAM converges with rate $r$ when there is no noise, we provide a feasible solution to the LMI~\eqref{eq:lmi_cvx_rate} with $\ell=1$. Our solution has the same structure as the weighted off-by-one IQC formulation in~\eqref{eq:IQC_parameterization}, where the positive definite matrix $Q\succ 0$ is given by
\begin{align*}
    Q &= \frac{m}{(3-r)(1-r)^2(1+r)^3(m-L\,(1-r)^2)}
      \bmat{q_{11} & q_{12} & q_{13} \\ q_{12} & q_{22} & q_{23} \\ q_{13} & q_{23} & q_{33}},\\
\text{with:}\quad    q_{11} &= r\,\bigl(2m^2(1-r)+2Lm\,(4+r-2r^2+r^3)-L^2(1-r^2)^2\bigr),\\
    q_{12} &= r\,\bigl(-2m^2(1-r)-2Lm\,(1+r)^2+L^2(4-r)(1-r^2)^2\bigr), \\
    q_{13} &= (3-r)(1-r^2)\bigl(-m\,(1+r^2)+L\,(1+r-2r^2-r^3+r^4)\bigr), \\
    q_{22} &= r\,\bigl(2m^2(1-r)-2Lm\,(2-3r-4r^2+r^3)+L^2(2-4r+r^2)(1-r^2)^2\bigr), \\
    q_{23} &= r\,(3-r)(1-r^2)\bigl(m\,(-1+2r+r^2)-L\,(-1+r-r^3+r^4)\bigr), \\
    q_{33} &= r\,(3-r)^2(1-r^2)^2.
  \end{align*}
  Using Mathematica~\cite{mathematica}, we verify that, for all $0<m\leq L$ with $1-\sqrt{\frac{m}{L}} \leq r < 1$, this is a feasible solution to a modified version of the LMI~\eqref{eq:lmi_cvx_rate} in which the term $\X_r^\tp \X_r$ is replaced by $\X_r^\tp T^\tp Q T \X_r$, where $T = \sbmat{A & B \\ -LC & 1}$. Since ${T^\tp Q T\succ 0}$, we have $\X_r^\tp T^\tp Q T \X_r \succeq c\,\X_r^\tp \X_r$, where $c > 0$ is the minimum eigenvalue of $T^\tp Q T$. Therefore, scaling the solution by $1/c$ yields a feasible solution to the original LMI~\eqref{eq:lmi_cvx_rate}. \cref{thm:lmi_cvx} then implies that RAM has convergence rate at most $r$. Since $\rho(A+mBC) = r$, the bound is achieved by $f(y) = \tfrac{m}{2} \|y\|^2$, so the convergence rate is \textit{exactly}~$r$.\qed

\subsection{Proof of \texorpdfstring{\cref{cor:param_restriction} (parameter restriction)}{parameter restriction}}\label{app:proof_param_bounds}

  The Jury criterion for stability ($r < 1$) is given in~\eqref{eq:jury_one}. Combining $\eqref{eb}+2\cdot\eqref{ed}$  with \eqref{ea}, we obtain: $0 < \alpha q  < 4$. For this to hold for all $q  \in [m,L]$, we must have $0 < \alpha < \frac{4}{L}$.
  Combining \eqref{ec} and \eqref{ed}, we obtain $-1+\alpha\eta q  < \beta < 1+\alpha\eta q$, which holds for all $q\in [m,L]$. When $\alpha\eta\geq0$, the $\beta$ range reduces to $-1+L(\alpha\eta) <\beta < 1+m(\alpha\eta)$ and thus, $\alpha\eta < \frac{2}{L-m}$. When $\alpha\eta< 0$, we instead obtain $-1+m(\alpha\eta) < \beta < 1+L(\alpha\eta)$ and $\frac{-2}{L-m} < \alpha\eta$.
  Although these bounds are derived for $Q_{m,L}$, the nestedness property $Q_{m,L} \subseteq F_{m,L}$ implies that these necessary conditions on $(\alpha,\beta,\eta)$ also hold for $F_{m,L}$.\qed

\newpage
\bibliographystyle{abbrv}
\bibliography{optalg}

\begin{thebibliography}{10}

\bibitem{katyusha}
Z.~Allen-Zhu.
\newblock Katyusha: {T}he first direct acceleration of stochastic gradient
  methods.
\newblock In {\em ACM SIGACT Symposium on Theory of Computing}, pages
  1200--1205, 2017.

\bibitem{antsaklis}
P.~J. Antsaklis and A.~N. Michel.
\newblock {\em Linear systems}.
\newblock Springer Science \& Business Media, 2006.

\bibitem{mosek}
M.~ApS.
\newblock {\em The MOSEK optimization suite 9.2.49}, 2021.

\bibitem{M-ASG}
N.~S. Aybat, A.~Fallah, M.~G{\"u}rb{\"u}zbalaban, and A.~Ozdaglar.
\newblock A universally optimal multistage accelerated stochastic gradient
  method.
\newblock In {\em Advances in Neural Information Processing Systems},
  volume~32, 2019.

\bibitem{mert}
N.~S. Aybat, A.~Fallah, M.~Gurbuzbalaban, and A.~Ozdaglar.
\newblock Robust accelerated gradient methods for smooth strongly convex
  functions.
\newblock {\em SIAM Journal on Optimization}, 30(1):717--751, 2020.

\bibitem{bassily2014}
R.~Bassily, A.~Smith, and A.~Thakurta.
\newblock Private empirical risk minimization: {E}fficient algorithms and tight
  error bounds.
\newblock {\em IEEE 55th Annu. Symp. on Found. of Computer Science}, 2014.

\bibitem{julia}
J.~Bezanson, A.~Edelman, S.~Karpinski, and V.~B. Shah.
\newblock Julia: {A} fresh approach to numerical computing.
\newblock {\em SIAM Review}, 59(1):65--98, 2017.

\bibitem{birand2013}
B.~Birand, H.~Wang, K.~Bergman, and G.~Zussman.
\newblock Measurements-based power control - a cross-layered framework.
\newblock In {\em Optical Fiber Communication Conference/National Fiber Optic
  Engineers Conference 2013}, 2013.

\bibitem{bottou2018}
L.~Bottou, F.~E. Curtis, and J.~Nocedal.
\newblock Optimization methods for large-scale machine learning.
\newblock {\em SIAM Review}, 60(2):223--311, 2018.

\bibitem{boyd_admm}
S.~Boyd, N.~Parikh, and E.~Chu.
\newblock {\em Distributed optimization and statistical learning via the
  alternating direction method of multipliers}.
\newblock Now Publishers Inc, 2011.

\bibitem{boyd_vandenberghe}
S.~Boyd and L.~Vandenberghe.
\newblock {\em Convex optimization}.
\newblock Cambridge University Press, 2004.

\bibitem{toulis_convergence_diagnostic}
J.~Chee and P.~Toulis.
\newblock Convergence diagnostics for stochastic gradient descent with constant
  learning rate.
\newblock In {\em International Conference on Artificial Intelligence and
  Statistics}, pages 1476--1485, 2018.

\bibitem{hypatia}
C.~Coey, L.~Kapelevich, and J.~P. Vielma.
\newblock Performance enhancements for a generic conic interior point
  algorithm.
\newblock {\em Mathematical Programming Computation}, 15:53--101, 2023.

\bibitem{cohen2018}
M.~Cohen, J.~Diakonikolas, and L.~Orecchia.
\newblock On acceleration with noise-corrupted gradients.
\newblock In {\em International Conference on Machine Learning}, pages
  1019--1028, 2018.

\bibitem{rmm}
S.~Cyrus, B.~Hu, B.~Van~Scoy, and L.~Lessard.
\newblock {A robust accelerated optimization algorithm for strongly convex
  functions}.
\newblock In {\em American Control Conference}, pages 1376--1381, 2018.

\bibitem{deklerk2017}
E.~De~Klerk, F.~Glineur, and A.~B. Taylor.
\newblock On the worst-case complexity of the gradient method with exact line
  search for smooth strongly convex functions.
\newblock {\em Optimization Letters}, 11:1185--1199, 2017.

\bibitem{deKlerk2020}
E.~De~Klerk, F.~Glineur, and A.~B. Taylor.
\newblock Worst-case convergence analysis of inexact gradient and newton
  methods through semidefinite programming performance estimation.
\newblock {\em SIAM Journal on Optimization}, 30(3):2053--2082, 2020.

\bibitem{saga}
A.~Defazio, F.~Bach, and S.~Lacoste-Julien.
\newblock {SAGA}: A fast incremental gradient method with support for
  non-strongly convex composite objectives.
\newblock In {\em Advances in Neural Information Processing Systems},
  volume~27, 2014.

\bibitem{devolder-strongly-convex}
O.~Devolder, F.~Glineur, and Y.~Nesterov.
\newblock First-order methods with inexact oracle: the strongly convex case.
\newblock Core discussion papers; 2013/16, Universit\'{e} Catholique de
  Louvain, 2013.

\bibitem{IGM}
O.~Devolder, F.~Glineur, and Y.~Nesterov.
\newblock Intermediate gradient methods for smooth convex problems with inexact
  oracle.
\newblock Core discussion papers; 2013/17, Universit\'{e} Catholique de
  Louvain, 2013.

\bibitem{dgn}
O.~Devolder, F.~Glineur, and Y.~Nesterov.
\newblock First-order methods of smooth convex optimization with inexact
  oracle.
\newblock {\em Mathematical Programming}, 146:37--75, 2014.

\bibitem{oracle-complexity}
Y.~Drori and A.~B. Taylor.
\newblock On the oracle complexity of smooth strongly convex minimization.
\newblock {\em Journal of Complexity}, 68(C), 2022.

\bibitem{drori-teboulle}
Y.~Drori and M.~Teboulle.
\newblock Performance of first-order methods for smooth convex minimization: a
  novel approach.
\newblock {\em Mathematical Programming}, 145:451--482, 2014.

\bibitem{fadali2013digital}
M.~S. Fadali and A.~Visioli.
\newblock {\em Digital control engineering: analysis and design}.
\newblock Academic Press, 2013.

\bibitem{spider}
C.~Fang, C.~J. Li, Z.~Lin, and T.~Zhang.
\newblock {SPIDER}: {N}ear-optimal non-convex optimization via stochastic path
  integrated differential estimator.
\newblock In {\em Advances in Neural Information Processing Systems}, pages
  687--697, 2018.

\bibitem{fazlyab}
M.~Fazlyab, A.~Ribeiro, M.~Morari, and V.~M. Preciado.
\newblock Analysis of optimization algorithms via integral quadratic
  constraints: {N}onstrongly convex problems.
\newblock {\em SIAM Journal on Optimization}, 28(3):2654--2689, 2018.

\bibitem{kakade_step_decay}
R.~Ge, S.~M. Kakade, R.~Kidambi, and P.~Netrapalli.
\newblock The step decay schedule: A near optimal, geometrically decaying
  learning rate procedure for least squares.
\newblock In {\em Advances in Neural Information Processing Systems},
  volume~32, 2019.

\bibitem{HB:convex}
E.~Ghadimi, H.~R. Feyzmahdavian, and M.~Johansson.
\newblock Global convergence of the {H}eavy-ball method for convex
  optimization.
\newblock In {\em 2015 European Control Conference}, pages 310--315, 2015.

\bibitem{ghadimi-lan1}
S.~Ghadimi and G.~Lan.
\newblock Optimal stochastic approximation algorithms for strongly convex
  stochastic composite optimization {I}: {A} generic algorithmic framework.
\newblock {\em SIAM Journal on Optimization}, 22(4):1469--1492, 2012.

\bibitem{ghadimi-lan2}
S.~Ghadimi and G.~Lan.
\newblock Optimal stochastic approximation algorithms for strongly convex
  stochastic composite optimization, {II}: {S}hrinking procedures and optimal
  algorithms.
\newblock {\em SIAM Journal on Optimization}, 23(4):2061--2089, 2013.

\bibitem{diss_ICML}
B.~Hu and L.~Lessard.
\newblock {Dissipativity theory for Nesterov's accelerated method}.
\newblock In {\em International Conference on Machine Learning}, pages
  1549--1557, 2017.

\bibitem{sgdn}
B.~Hu, P.~Seiler, and L.~Lessard.
\newblock Analysis of biased stochastic gradient descent using sequential
  semidefinite programs.
\newblock {\em Mathematical Programming}, 187:383--408, 2020.

\bibitem{kakade_accelerated_stochastic}
P.~Jain, S.~M. Kakade, R.~Kidambi, P.~Netrapalli, and A.~Sidford.
\newblock Accelerating stochastic gradient descent for least squares
  regression.
\newblock In {\em Conference On Learning Theory}, volume~75, pages 545--604,
  2018.

\bibitem{svrg}
R.~Johnson and T.~Zhang.
\newblock Accelerating stochastic gradient descent using predictive variance
  reduction.
\newblock In {\em Advances in Neural Information Processing Systems},
  volume~26, 2013.

\bibitem{kulunchakov2019}
A.~Kulunchakov and J.~Mairal.
\newblock A generic acceleration framework for stochastic composite
  optimization.
\newblock In {\em Advances in Neural Information Processing Systems},
  volume~32, 2019.

\bibitem{lan2012}
G.~Lan.
\newblock An optimal method for stochastic composite optimization.
\newblock {\em Mathematical Programming}, 133:365--397, 2012.

\bibitem{lessardCSM}
L.~Lessard.
\newblock The analysis of optimization algorithms: A dissipativity approach.
\newblock {\em IEEE Control Systems Magazine}, 42(3):58--72, 2022.

\bibitem{lessard16}
L.~Lessard, B.~Recht, and A.~Packard.
\newblock Analysis and design of optimization algorithms via integral quadratic
  constraints.
\newblock {\em SIAM Journal on Optimization}, 26(1):57--95, 2016.

\bibitem{algosyn_ACC}
L.~Lessard and P.~Seiler.
\newblock {Direct synthesis of iterative algorithms with bounds on achievable
  worst-case convergence rate}.
\newblock In {\em American Control Conference}, 2020.

\bibitem{JUMP}
M.~Lubin, O.~Dowson, J.~{Dias Garcia}, J.~Huchette, B.~Legat, and J.~P. Vielma.
\newblock {JuMP} 1.0: {R}ecent improvements to a modeling language for
  mathematical optimization.
\newblock {\em Mathematical Programming Computation}, 2023.

\bibitem{iqc}
A.~Megretski and A.~Rantzer.
\newblock System analysis via integral quadratic constraints.
\newblock {\em IEEE Transactions on Automatic Control}, 42(6):819--830, 1997.

\bibitem{scherer}
S.~Michalowsky, C.~Scherer, and C.~Ebenbauer.
\newblock Robust and structure exploiting optimisation algorithms: an integral
  quadratic constraint approach.
\newblock {\em International Journal of Control}, 94(11):2956--2979, 2021.

\bibitem{mihailo}
H.~Mohammadi, M.~Razaviyayn, and M.~R. Jovanovi{\'c}.
\newblock Robustness of accelerated first-order algorithms for strongly convex
  optimization problems.
\newblock {\em IEEE Transactions on Automatic Control}, 66(6):2480--2495, 2021.

\bibitem{nesterov_book}
Y.~Nesterov.
\newblock {\em Lectures on convex optimization, second edition}, volume 137.
\newblock Springer, 2018.

\bibitem{sarah}
L.~M. Nguyen, J.~Liu, K.~Scheinberg, and M.~Tak\'{a}\v{c}.
\newblock {SARAH}: {A} novel method for machine learning problems using
  stochastic recursive gradient.
\newblock In {\em International Conference on Machine Learning}, pages
  2613--2621, 2017.

\bibitem{nocedal_wright}
J.~Nocedal and S.~Wright.
\newblock {\em Numerical optimization}.
\newblock Springer Science \& Business Media, 2006.

\bibitem{polyak_book}
B.~T. Polyak.
\newblock {\em Introduction to optimization}.
\newblock Translations series in mathematics and engineering. Optimization
  Software, Inc., 1987.

\bibitem{ryu2020}
E.~Ryu, A.~B. Taylor, C.~Bergeling, and P.~Giselsson.
\newblock Operator splitting performance estimation: {T}ight contraction
  factors and optimal parameter selection.
\newblock {\em SIAM Journal on Optimization}, 30(3):2251--2271, 2020.

\bibitem{scherer2021convex}
C.~Scherer and C.~Ebenbauer.
\newblock Convex synthesis of accelerated gradient algorithms.
\newblock {\em SIAM Journal on Control and Optimization}, 59(6):4615--4645,
  2021.

\bibitem{taylor2019}
A.~B. Taylor and F.~Bach.
\newblock Stochastic first-order methods: non-asymptotic and computer-aided
  analyses via potential functions.
\newblock In {\em Conference on Learning Theory}, volume~99, pages 2934--2992,
  2019.

\bibitem{ITEM}
A.~B. Taylor and Y.~Drori.
\newblock An optimal gradient method for smooth strongly convex minimization.
\newblock {\em Mathematical Programming}, 199:557--594, 2022.

\bibitem{taylor2017smooth}
A.~B. Taylor, J.~M. Hendrickx, and F.~Glineur.
\newblock Smooth strongly convex interpolation and exact worst-case performance
  of first-order methods.
\newblock {\em Mathematical Programming}, 161:307--345, 2017.

\bibitem{taylor2018exact}
A.~B. Taylor, J.~M. Hendrickx, and F.~Glineur.
\newblock Exact worst-case convergence rates of the proximal gradient method
  for composite convex minimization.
\newblock {\em Journal of Optimization Theory and Applications},
  178(2):455--476, 2018.

\bibitem{taylor2018lyapunov}
A.~B. Taylor, B.~Van~Scoy, and L.~Lessard.
\newblock Lyapunov functions for first-order methods: {T}ight automated
  convergence guarantees.
\newblock In {\em International Conference on Machine Learning}, pages
  4897--4906, 2018.

\bibitem{tmm}
B.~Van~Scoy, R.~A. Freeman, and K.~M. Lynch.
\newblock The fastest known globally convergent first-order method for
  minimizing strongly convex functions.
\newblock {\em IEEE Control System Letters}, 2(1):49--54, 2017.

\bibitem{alift_cdc}
B.~Van~Scoy and L.~Lessard.
\newblock {Absolute stability via lifting and interpolation}.
\newblock In {\em IEEE Conference on Decision and Control}, pages 6217--6223,
  2022.

\bibitem{tut-lyaplift_cdc}
B.~Van~Scoy and L.~Lessard.
\newblock {A tutorial on a Lyapunov-based approach to the analysis of iterative
  optimization algorithms}.
\newblock In {\em IEEE Conference on Decision and Control}, pages 3003--3008,
  2023.

\bibitem{wang2013}
C.~Wang, X.~Chen, A.~J. Smola, and E.~P. Xing.
\newblock Variance reduction for stochastic gradient optimization.
\newblock In {\em Advances in Neural Information Processing Systems},
  volume~26, 2013.

\bibitem{mathematica}
{Wolfram Research, Inc.}
\newblock Mathematica, {V}ersion 12.3.1.
\newblock Champaign, IL, 2021.

\bibitem{zames-falb}
G.~Zames and P.~Falb.
\newblock Stability conditions for systems with monotone and slope-restricted
  nonlinearities.
\newblock {\em SIAM Journal on Control}, 6(1):89--108, 1968.

\bibitem{snvrg}
D.~Zhou, P.~Xu, and Q.~Gu.
\newblock Stochastic nested variance reduction for nonconvex optimization.
\newblock {\em Journal of Machine Learning Research}, 21(103):1--63, 2020.

\end{thebibliography}

\end{document}